\documentclass[twoside,leqno]{article}

\usepackage[letterpaper]{geometry}

\usepackage{AKstyle}

\usepackage{siamproceedings}

\usepackage[T1]{fontenc}
\usepackage{amsfonts}
\usepackage{graphicx}
\usepackage{epstopdf}
\usepackage{enumitem}
\usepackage{algorithmic}
\ifpdf
  \DeclareGraphicsExtensions{.eps,.pdf,.png,.jpg}
\else
  \DeclareGraphicsExtensions{.eps}
\fi

\usepackage{amsopn}

\begin{document}

\newcommand\relatedversion{}
\title{\Large The Shape of Math To Come\footnote{The title is an homage to Ornette Coleman's 1959 album, ``The Shape of Jazz to Come''.}}
    \author{Alex Kontorovich\thanks{Department of Mathematics, Rutgers University (\email{alex.kontorovich@rutgers.edu}, \url{http://math.rutgers.edu/~alexk}). Partially supported by NSF grant DMS-2302641, US-Israel BSF grant 2020119, and a Simons Foundation Fellowship.}
    }

\date{}

\maketitle

\fancyfoot[R]{\scriptsize{Copyright \textcopyright\ 2025 by SIAM\\
Unauthorized reproduction of this article is prohibited}}

\fancyfoot[C]{\thepage}

\begin{abstract}
We present an overview of how certain computational tools currently interact with mathematical practice, and reflect on the implications for research mathematics in the short to medium term, as the field navigates the emerging age of AI and formal verification systems.
\end{abstract}

\section{Introduction.}\

This paper will discuss a vision for what research mathematics may look like in the age we now seem to be entering, of AI and formalization. Given the pace of progress, my words will likely be obsolete almost the instant they appear, but I will nevertheless take the risk of setting down how I see the landscape, in September 2025.
\\

My own involvement in formalization began about five years ago (see \cite{GomesKontorovich2020}), driven by the particular challenges of my research area. Analytic number theory is characterized by particularly long, technical arguments spanning scores of pages. These proofs typically synthesize numerous lemmas, each valid only within constrained ranges of their parameters, which must align perfectly at the end for the main theorem to hold. This intricate dance of parameters and constraints is both what I love about the subject, and what makes it treacherous: a single minus sign error buried somewhere deep in the middle of one of these technical arguments can invalidate the entire main result. 

When an argument involves something like a complicated algebraic manipulation, 
I can use a computer algebra system such as Sage or Mathematica to verify it once and for all; then I no longer waste time going over that part of the argument again and again to ensure no errors were made. Shouldn't there be similar software, but to track the logical dependencies and parameter constraints across my entire paper? 

Such software does exist, and is called an Interactive Theorem Prover; see \secref{sec:Lean}. It, in principle, allows users to ``formalize'' their papers, entering the proof line by line until the computer certifies that the argument is complete. The reality, however, is sobering. As things currently stand, I cannot even \emph{state} formally most theorems of interest to me, much less prove them. The existing formal libraries (see \secref{sec:Mathlib}) would need to grow by orders of magnitude before doing research on top of them becomes a reality.

This leads to a fundamental conundrum: new mathematics seems to grow at an exponential rate, say, $\delta$, and formal libraries also grow exponentially, at a rate, say, $\varepsilon$; but as far as I can tell, we currently have:
$$
\delta > \varepsilon.
$$
So if this inequality persists, formalization will \emph{never} catch up to natural language mathematics research, and the dream of software as useful to \emph{reasoning} as computer algebra systems are to symbolic manipulation will forever remain a distant mirage.

Unless, that is, formalization can be ``supercharged'' by automation. And so I've gone from wanting to do research math, to wanting sufficient formalization to do research math, to wanting sufficient AI to do sufficient formalization to do research math!  What follows is my vision for how we might get there, and how things might look if/when we do.
The old adage: ``It’s hard to make predictions, especially about the future,'' variously attributed to Yogi Berra, Niels Bohr, and others, is very apt here. So before we begin, let's take a step back.

\section{In the Year 2000...}\

At the turn of the millennium, some of the world’s most prominent mathematicians gathered in Tel Aviv for a conference to discuss what mathematics might look like in the 21st century. The proceedings were published as a special issue in GAFA called ``Visions of Mathematics'', and I strongly recommend reading the whole issue \cite{GAFA2000}. A particularly prescient contribution was given by Tim Gowers, titled ``Rough Structure and Classification'' \cite{Gowers2000}.

Section 2 of Gowers’s paper is provocatively called: ``Will Mathematics Exist in 2099?'' In it, he imagines a ``dialogue between a mathematician and a computer in two to three decades’ time''  --  and here we are! His imagined dialogue is eerily similar to the experiences so many of us  now have on a regular basis when interacting with ``chatbots'' such as OpenAI's GPT, Google's Gemini, Anthropic's Claude, xAI's Grok, and others:

\begin{figure}[h]
    \centering
\includegraphics[width=.7\textwidth]{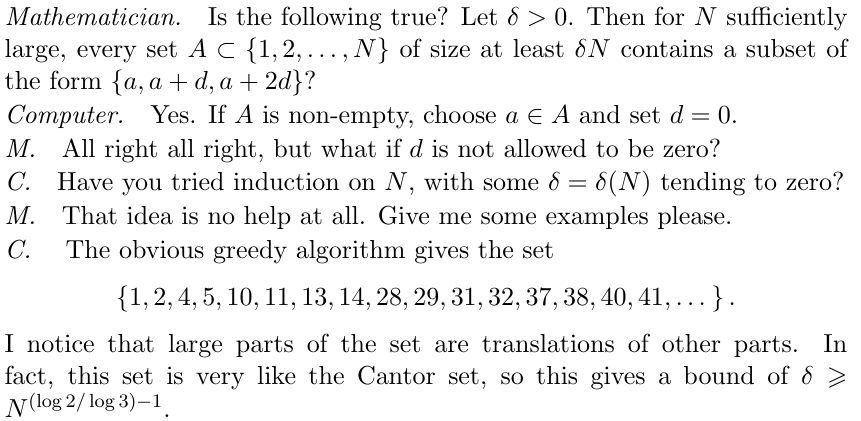}
    \caption{A snippet from \cite{Gowers2000}}
    \label{fig:gowersSnippet}
\end{figure}

Gowers notes: ``The idea of the dialogue is that the computer is very helpful to the mathematician, while not doing anything \emph{particularly clever}'' (emphasis added). So, in that vein: let's establish some common ground and be precise about what I shall mean by a ``Large Language Model'' (LLM).

\section{What is an LLM?}\

We've all had the following experience: when reading a book, you reach the bottom of a page mid-sentence, and as your fingers fumble, struggling to turn to the next page, your mind can't help but race ahead, filling in how the sentence might continue. For example:

\begin{figure}[h]
    \centering
\includegraphics[width=.55\textwidth]{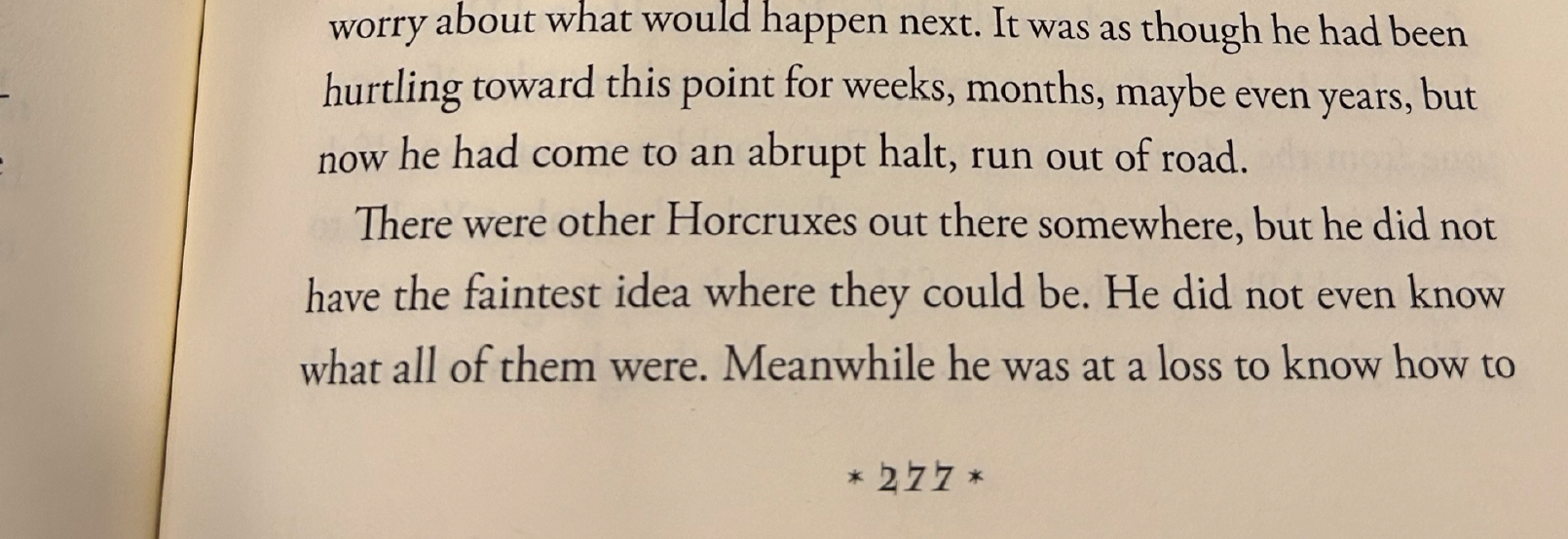}
    \caption{The bottom of a page from ``Harry Potter and the Deathly Hallows'' by J. K. Rowling}
    \label{fig:HP}
\end{figure}

We read that
\begin{quote}
``[Harry] did not even know what all [the Horcruxes] were. Meanwhile he was at a loss to know how to\ldots''
\end{quote}
How to \emph{what}? Your mind instinctively begins to generate possibilities. How to ``destroy'' the Horcruxes? Or perhaps how to ``find'' them? Or, far less plausibly, how to ``sandwich''?

I imagine standing at a game of Family Feud, hearing Steve Harvey announce: ``Survey says!\ldots'' But here the ``survey'' doesn't poll 100 random people -- it samples all available human text! The results might look something like this: given the context, roughly 20\% of people, say, would continue the sentence with ``destroy,'' another 7\% would use, ``find,'' and perhaps a tiny 0.003\% would suggest ``sandwich.''

\emph{This}, at its core, is a Large Language Model: it is a function whose input is ``context'' -- a possibly very long sequence of ``tokens'' (words, punctuation marks, etc.) -- and whose output is \emph{not} the next word (as is commonly misconstrued), but rather a \emph{probability distribution} over all possible next tokens. 

\begin{figure}[h]
    \centering
\includegraphics[width=.22\textwidth]{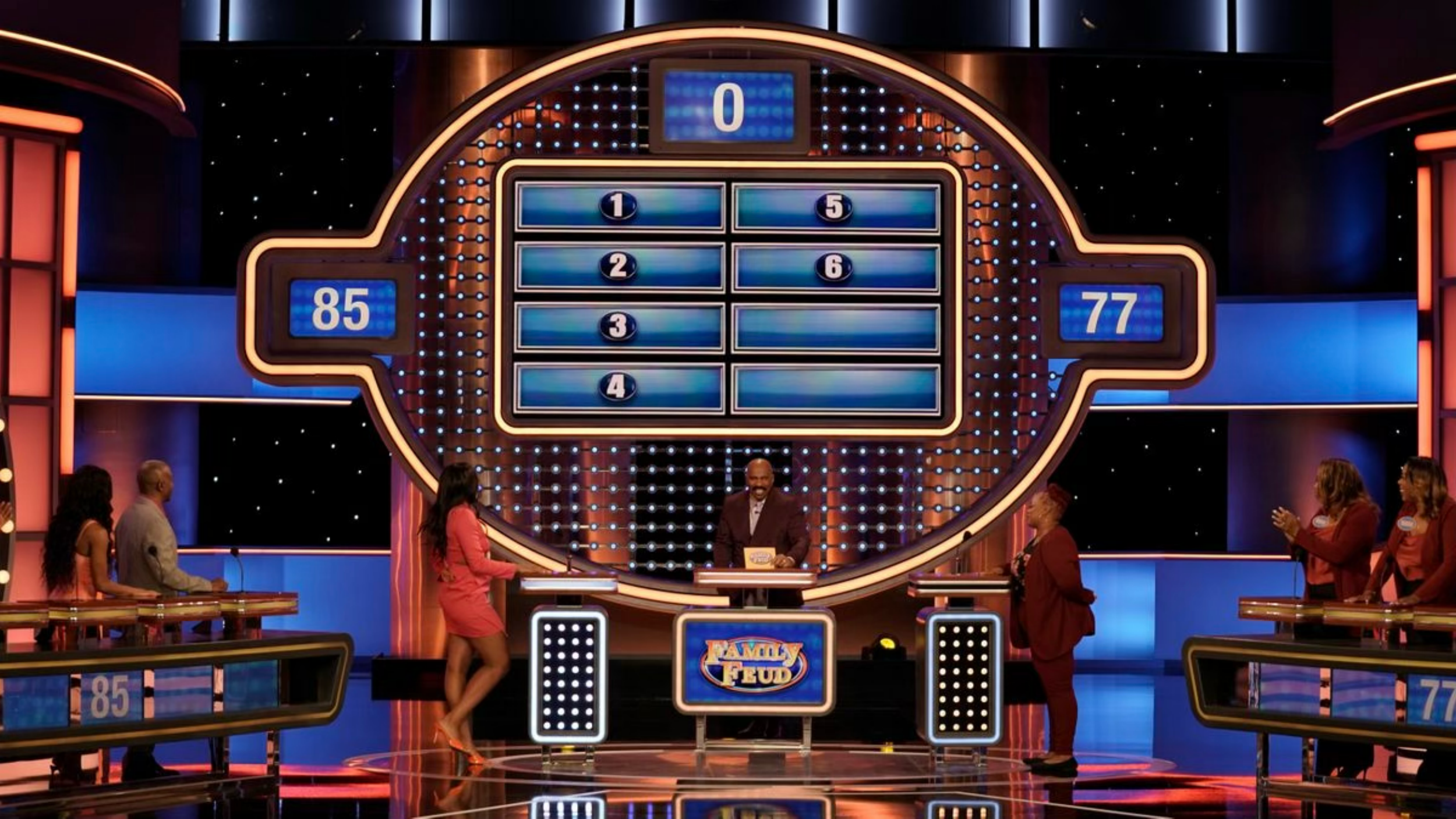} \hskip.1in
\includegraphics[width=.75\textwidth]{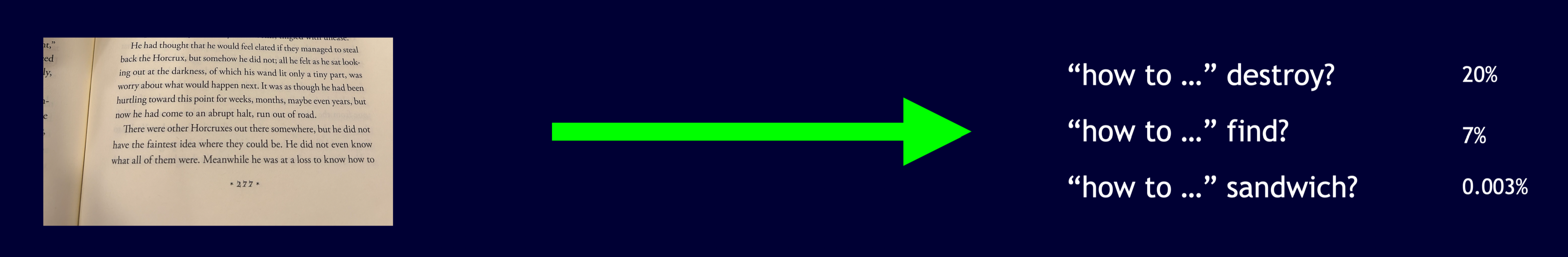}
    \caption{An LLM is a function from context (a sequence of tokens) to a likelihood distribution of potential next tokens}
    \label{fig:LLM}
\end{figure}

For our purposes, we can set aside the unfathomable engineering that makes such a function possible: the architecture (transformers), the representation of data (embeddings in high-dimensional vector spaces), the training process (gradient descent, pre-training, fine-tuning), and so on. We simply take as given that such a function exists.

It is on top of such a model that one builds the aforementioned chatbot; the bot must decide how to \emph{choose} the actual next token from the probability distribution the model outputs. This choice depends on parameters like the ``temperature'' -- perhaps always selecting the most likely token (which tends to produce somewhat stilted prose), or sampling randomly from the distribution, or various options in between. Once a token is chosen, the bot appends it to the context, and feeds the entire sequence back through the model; it then chooses the next token, and the next and the next, until it selects a special token that signals ``stop talking.''
(We will return later to more complex agentic behavior, like use of tools and self-directed workflows; for now, I’m focusing on the mechanism of a ``pure'' chatbot.) It is perhaps this pure chatbot behavior that Gowers was envisioning as a computer process ``not doing anything particularly clever'' in its discussion with you. 

A crucial observation here is that this process is profoundly \emph{stochastic}. At each step, the bot makes a random choice, and each choice affects the probability distribution for all subsequent choices. The randomness compounds with every token -- random upon random upon random. This stochasticity is simultaneously the system's great strength -- enabling it to tackle a vast range of problems with remarkable versatility -- but as we shall see, it creates fundamental challenges when applied to deterministic problems like mathematical proof.

\section{Stochastic vs. Deterministic Problems.}\

Consider the difference between approximate and exact correctness. A long poem with one poorly-chosen word in the middle can still be beautifully moving; the same for an improvised jazz solo with an accidentally dissonant note. These are domains where 99\% success represents genuine achievement. But a 99\% correctly proved theorem? A single error can render everything that follows meaningless! Mathematics requires deterministic correctness, not probabilistic approximation.

To understand why this matters quantitatively, suppose an LLM achieves 99\% accuracy at each individual step -- a remarkably high success rate. If a mathematical argument requires chaining 1000 such steps in sequence, with the success of each step independent of the others, what is the probability of overall success? The answer: $0.99^{1000} \approx 0.004\%$. A process that seems highly reliable at each step becomes almost certainly wrong when many steps compound.

The human brain apparently also operates in two rather distinct modes, popularized in Daniel Kahneman's \emph{Thinking, Fast and Slow} as ``System I'' versus ``System II''. The former is fast, automatic, and effortless; the latter is slow, deliberate, and effortful. When asked ``What is $8+3$?'' most adults respond instantly -- this is System I processing. You've encountered this calculation so many times that your brain has essentially memorized the answer. (This is why children learning their times tables is so important: it converts what was once effortful System II calculation into automatic System I recall, freeing cognitive resources for more complex tasks.)

But when I calculate $437 \times 82$ on paper, I require System II processing. However, because I've learned the standard algorithm, I can decompose this problem into a long chain of simple System I calculations: first compute $2 \times 7 = 14$, write the $4$, carry the $1$, and so on. The key insight is that System II mathematical reasoning often consists of carefully chaining together many System I steps, where success requires \emph{every link in the chain} to be correct.

This brings us back to our stochastic problem. If your ``System I machine'' is 99\% accurate at each step, but you need to chain 1000 steps to solve a System II problem, you face the aforementioned $0.004\%$ success rate. This is why using stochastic tools for deterministic mathematics is fundamentally challenging.

Yet this argument assumes that AI capabilities will remain static. Given the extraordinary pace of progress in recent years, perhaps these limitations are merely temporary?

\section{AI Capabilities Are Advancing Rapidly...}\

The evidence for continued rapid progress is compelling. Consider performance on the International Mathematical Olympiad (IMO), widely regarded as one of the most challenging mathematical competitions for high school students:

\begin{itemize}
\item In 2023, the number of AI systems that could earn a single point on the IMO was... zero.
\item In 2024, DeepMind's AlphaProof + AlphaGeometry achieved a silver medal performance, one point shy of gold \cite{AlphaProof2024}.
\item In 2025, multiple AI systems claimed gold medal scores, with roughly half using formal verification methods and half working purely in natural language.\footnote{While DeepMind used a general-purpose natural-language model, it still called an AlphaGeometry-type system for the geometry problems.}
\end{itemize}

This trajectory is remarkable. Within two years, we moved from complete inability to gold medal performance. Moreover, the 2025 success of purely natural language approaches demonstrates that formal verification is not strictly necessary for solving difficult mathematical problems -- at least at the IMO level.

By the next IMO in the summer of 2026, it may be reasonable to predict that the number of leading AI systems competing in the IMO will return to zero -- not because they've regressed, but because the competition will have become trivially easy for them. (Open source models may continue to participate, but commercial systems will likely find it beneath their capabilities.)

This pattern of accelerating capability is not unique to mathematics. It reflects what Richard Sutton famously called the ``Bitter Lesson'' of seventy years of AI research: domain-specific solutions crafted by human experts are consistently overtaken by general-purpose methods that leverage ever-increasing computational resources \cite{Sutton2019}. Time and again, we've seen that scaling compute and data beats clever engineering tricks.

If Moore's Law continues -- or even if it merely slows rather than stops -- we can expect LLMs to sustain progressively longer chains of reasoning. Systems that struggle with 100-step arguments today may handle 1000-step arguments tomorrow, and 10,000-step arguments the day after. The transition from high school mathematics to undergraduate to PhD-level to research-level mathematics may simply require scaling up the sustained reasoning capability by some large multiple.

Indeed, one can pose the following a loose analogy to chess ratings: a 2500-ELO rated chess grandmaster might correspond to a freshly minted mathematics PhD; 2600 to an established researcher; 2700 to a leading expert; and 2800+ to mathematical giants. If AI systems can reach the 2300-equivalent level for elite high school mathematics, why not eventually reach 2800 and beyond?

So perhaps the stochastic/deterministic tension I've emphasized will simply dissolve as systems become powerful enough to maintain reliability across arbitrarily long chains of reasoning. Perhaps natural language AI will, after all, be sufficient for all of mathematics?
It seems to me that even this optimistic scenario fails to solve a fundamental problem.

\section{... But Even Perfect Scaling Is Insufficient.}\

Let us imagine, optimistically, that AI continues to improve at its current pace. Suppose that within a decade, we have systems capable of producing novel mathematical research -- not just solving competition problems, but discovering new theorems, developing new theories, and writing papers that advance human knowledge.

Consider a future AI system that can generate 100 research papers per day. To make the scenario as favorable as possible, suppose that this system has achieved extraordinary reliability: 99 of those 100 papers are completely correct, with rigorous proofs and meaningful contributions to mathematics. Only 1 paper in 100 contains an error -- perhaps a subtle mistake in a lemma, or a case that wasn't properly handled, or a minus sign buried somewhere deep in a calculation.

Would you rely on such a system? I would \textbf{not}.

To use this tool effectively, I would need to read through all 100 papers to determine which 99 are correct and which 1 is flawed. But verifying the correctness of a mathematical paper may not necessarily be substantially easier than producing it in the first place, especially for cutting-edge research. There would certainly be value in mining the tool's output for interesting ideas, novel approaches, or unexpected connections; however, I could not treat the theorems themselves as established facts to build upon. I would need to spend my entire day, every day, checking these papers, only to discover that most of my effort went toward verifying work that was already correct, while the single flawed paper might pass undetected if its error is subtle enough.

One might object that this hypothetical 99\% accuracy is already better than the current published literature! After all, most mathematicians have encountered papers with errors, gaps in arguments, or cases where ``the details are left to the reader.'' This is all true. But there is a crucial difference: with human-written papers, we develop intuitions about reliability. As Kevin Buzzard wrote recently \cite{BuzzardPrivate2025}: ``LLMs would be very happy to paper over any or all cracks in the argument or even lie in order to get to the claims of the result and as such I feel far less motivated to read the paper properly on the basis that it was not 100\% generated by someone who believed that what they are writing is \emph{true}, but rather by something who believed that what they are writing is \emph{plausible}.'' (Emphasis added.)

We develop intuitions about reliability through multiple signals: the quality of exposition, whether techniques are applied in standard ways, how well edge cases are handled, whether claims seem plausible given known results, etc., etc. A brilliant proof from an unknown mathematician can be immediately recognized as such precisely because we can evaluate these qualities. (For just one such example, see \cite{Zhang2014a}.) With AI producing 100 papers daily from a single source, these evaluation heuristics provide no differentiation; every paper must be treated with equal suspicion.

Moreover, even for the 99 correct papers, I cannot extract their full value without verification. A beautiful theorem might suggest new directions for research, but I cannot confidently build on it until I've checked the proof myself. The ideas might be worth mining, the conjectures worth pursuing, but the theorems themselves remain uncertain until verified. And at that point, why not just work out the mathematics myself?

This is not a problem that scaling can solve. Even if we improve the accuracy to 99.99\% or 99.999\%, the fundamental issue remains: for mathematics, a non-zero error rate creates an unusable tool. If we cannot distinguish the correct results from the flawed ones without verification, then that defeats the whole purpose of automation.

Compare this to a different scenario: an AI that produces only 10 papers per year, but each comes with a \emph{formally verified proof}. Moreover -- and this is crucial -- suppose that I have learned to read formal statements well enough to verify that they say exactly what I want them to say, with no mistranslations or misinterpretations. \emph{That} is a tool I would use enthusiastically.

The difference is deterministic certification. I can trust those 10 papers completely, build upon them confidently, and incorporate their results into my work without fear of hidden errors. The slower pace is not merely acceptable -- it is a feature, not a bug. Quality and certainty matter far more than quantity when building a cumulative body of mathematical knowledge.

\section{Formally Verified Mathematics.}\label{sec:Lean}\

Interactive Theorem Provers (also called ``proof assistants'') are software systems that check every step of a mathematical proof against axioms and previously established results using the rules of logical inference. Such technology was envisioned in some form as early as Hilbert, if not Euclid, and has existed in practice since the late 1960s. Established systems include AUTOMATH, Mizar, Coq (now Rocq), Isabelle, HOL Light, and many others.

Among the most recent is Lean, which has been rapidly gaining adoption not just in computer science and mathematical logic circles, but in mainstream research mathematics. Many of the principles discussed here apply to any proof assistant, but I will focus on Lean and its ecosystem. (Full disclosure: I serve on the Strategic Advisory Board of the Lean Focused Research Organization.) 
One might worry that there are inherent limits to what can be expressed in Lean, but as Buzzard’s 2022 ICM Plenary Lecture \cite{Buzzard2023} argues, the formal framework has already proved itself broad enough to accommodate essentially any mathematics we may want to formalize; let me not repeat the evidence here.

There are some common misconceptions about what the term ``Lean'' actually refers to, so let me be precise.

\begin{figure}[ht]
    \centering
    \includegraphics[width=0.5\linewidth]{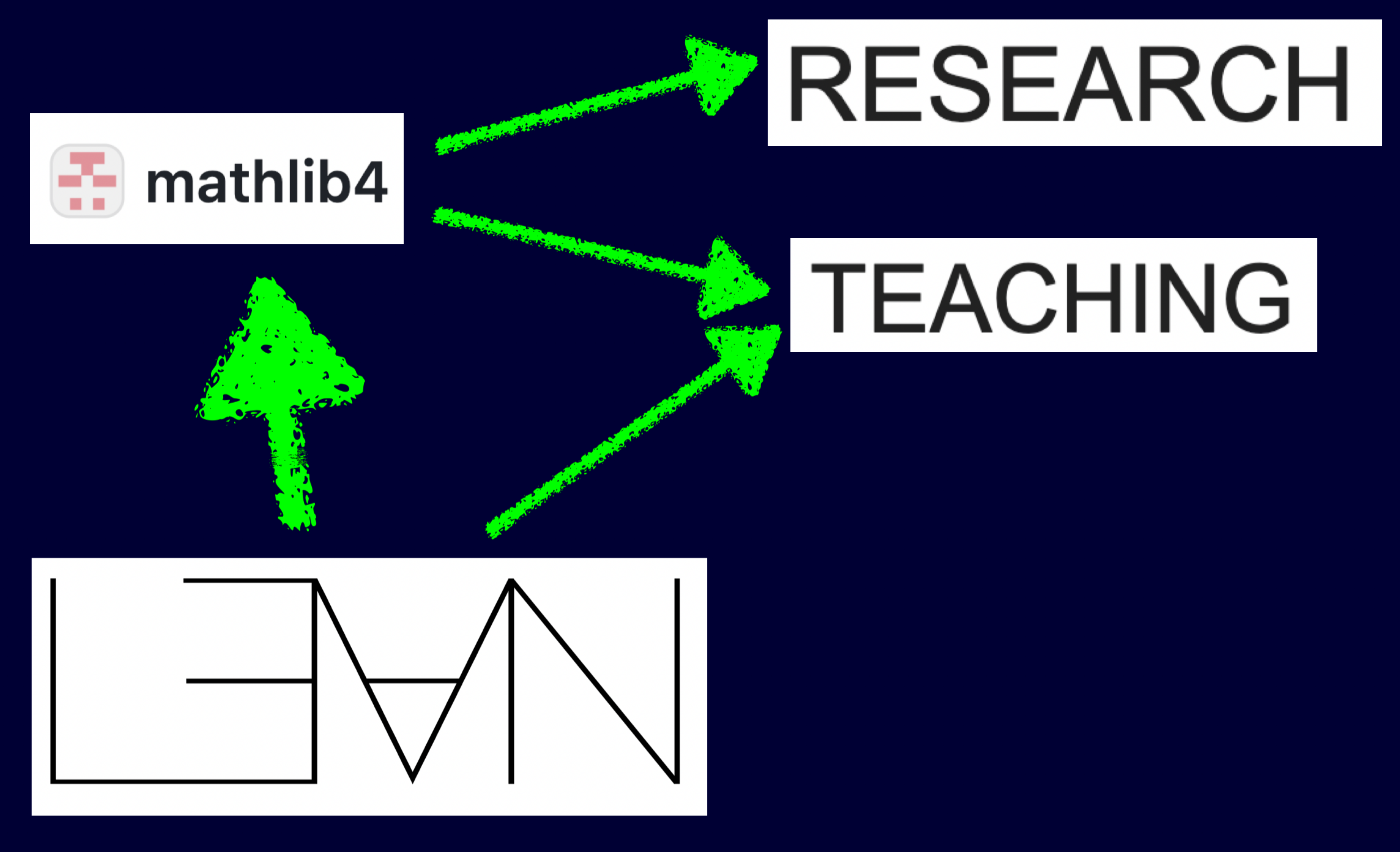}
    \caption{Different aspects of ``Lean''}
    \label{fig:Lean}
\end{figure}

\subsection{What is Lean?}

Core Lean (now technically Lean4) is software originally developed by Leonardo de Moura (then at Microsoft Research, now at Amazon Web Services) to certify bug-free code. Via the so-called Curry-Howard correspondence \cite{Curry1934, Howard1980}, this same technology can certify that a mathematical statement has been proved, checking every step down to axioms. It has nothing, \emph{a priori}, to do with AI.

The workflow is simple: humans write proofs in Lean's formal language, line by line. 
After each line, Lean either accepts the step as valid or indicates an error.
 When a proof is complete, Lean displays ``No goals'', meaning there are no remaining statements to prove.

Here is a simple example:

\begin{figure}[h]
    \centering
    
\raisebox{1.6in}{$(a)$}\hskip.1in\includegraphics[width=0.83\linewidth]{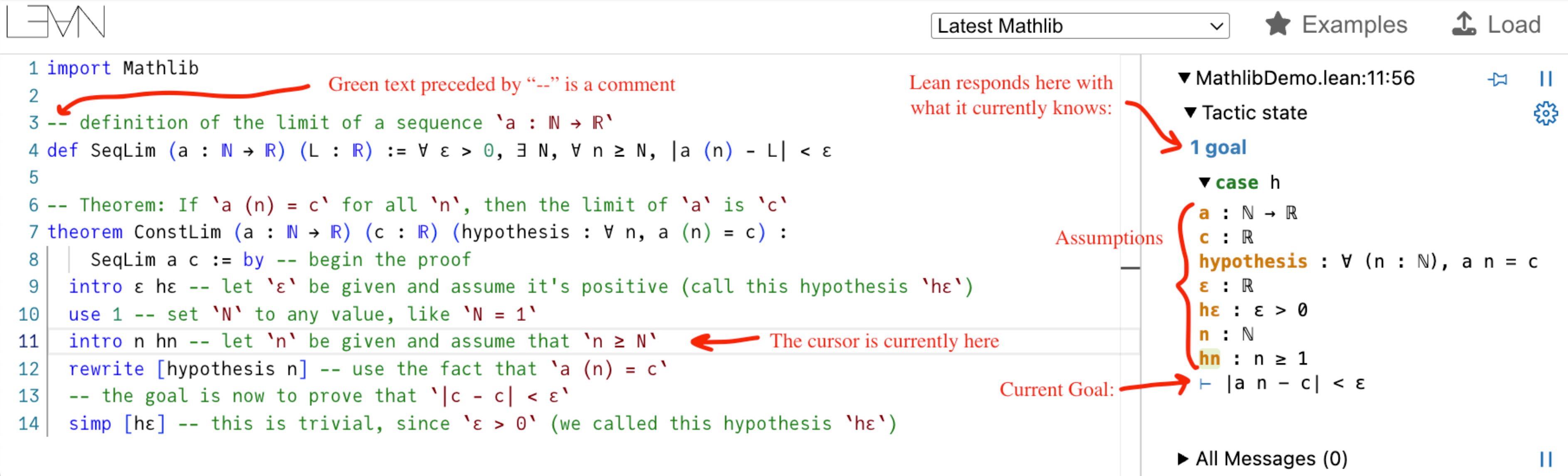}\\

\raisebox{1.6in}{$(b)$}\hskip.1in\includegraphics[width=0.83\linewidth]{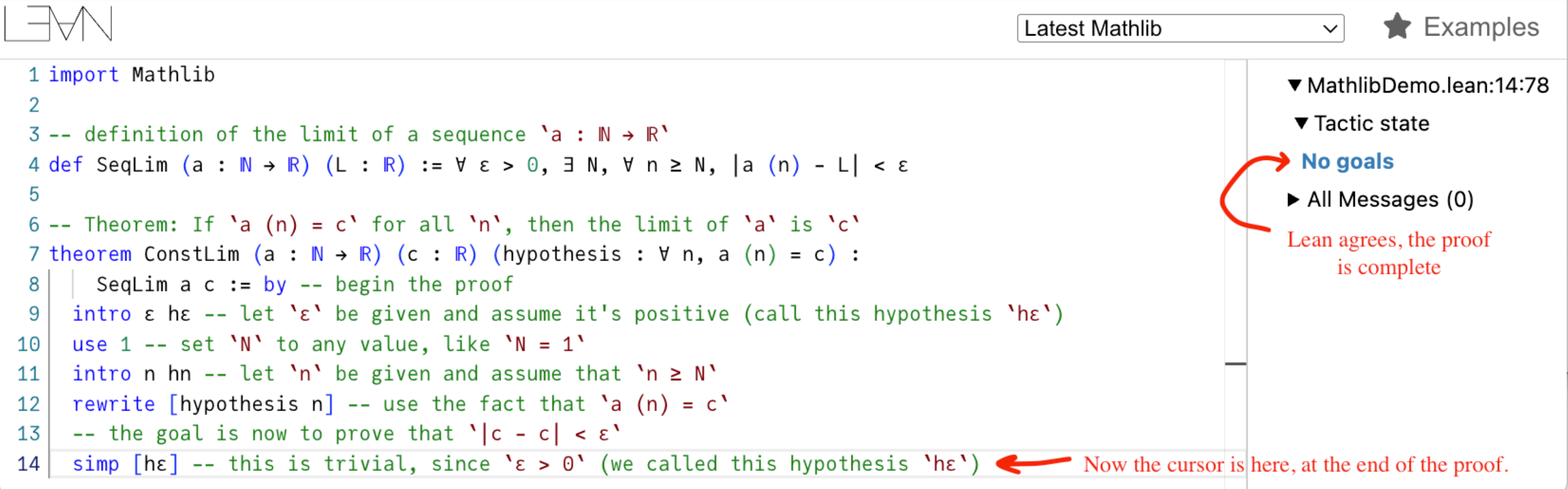}
    \caption{Proving a theorem in Lean.}
    \label{fig:Lean12}
\end{figure}

\subsection{What is Mathlib?}\label{sec:Mathlib}

Even if you have never tried to prove a significant theorem down to its axioms from scratch, you can surely imagine that it might be extraordinarily tedious and time consuming.
To make any progress, we need a comprehensive library of previously proved results, stated in maximal generality to be useful in the widest range of situations. We also need powerful ``tactics'' -- automated procedures that can solve entire classes of subgoals -- so that users can appeal to just a few such tactics in their proofs, rather than trying to memorize thousands of individual theorem names.

This is Mathlib: a massive, interconnected, and maintainable library built on Lean.
Like Lean itself, Mathlib is entirely open source, enabling the collaborative effort of hundreds of contributors worldwide.
When my mathematician friends complain that they looked at Mathlib and understood nothing despite being experts in the relevant field, I explain that Mathlib is \emph{not} meant as an educational project. Its single-minded goal is comprehensiveness and scalability.

\begin{figure}[h]
    \centering
    \includegraphics[width=0.7\linewidth]{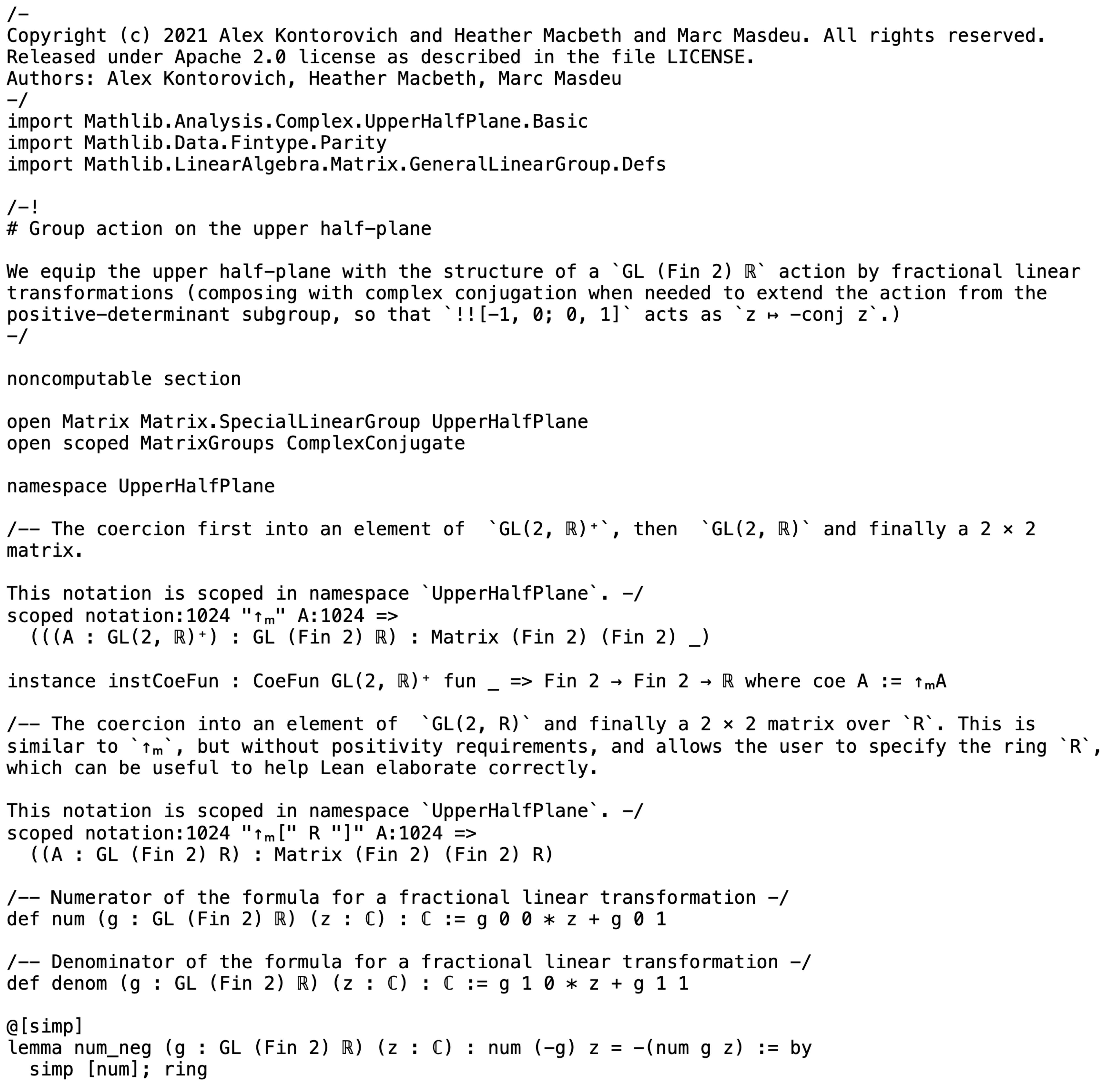}
    \caption{A snippet from a sample file \cite{KontorovichMacbethMasdeu2021} in Mathlib}
    \label{fig:Mathlib}
\end{figure}

Consider someone working on algebraic groups of Lie type. The Manifolds library must integrate seamlessly with the libraries on Groups and Algebraic Geometry, with no ``namespace'' conflicts or incompatible conventions. The Mathlib maintainers work extremely hard to ensure that you can import all these libraries simultaneously and get to work.

Mathlib must scale from its current order of magnitude (roughly 2 million lines of code) to 10 million, to 100 million lines and beyond, while maintaining interactive responsiveness. If Lean takes more than a few seconds to respond after you write a line of code, it becomes unusable for research mathematics.\\

It's worth emphasizing, as illustrated in \figref{fig:Lean}, that engaging with formal mathematics does not require contributing to Mathlib itself. If the mathematical foundations you need are already available in Mathlib, you can conduct research by building on top of the library without necessarily aiming to extend it. 
Similarly, Mathlib serves as a resource for teaching -- though for pedagogical purposes, it's often better to work directly from Lean's core library, as Mathlib's highly general and sophisticated treatments can be too advanced for students first learning the subject; see \secref{sec:Teaching}. 

We will return to both of these issues in more detail later. For now, suffice it to say that the three applications -- research using existing formalizations, teaching with appropriate levels of abstraction, and extending Mathlib with new mathematics -- are distinct activities that all serve different communities and purposes. Mathlib's role is to provide the foundational infrastructure that makes the other two possible.

\subsection{The Library Gap and AI Assistance.}

Now we can recall and better appreciate the fundamental conundrum from the introduction: research mathematics grows at rate $\delta$ while Mathlib grows at rate $\varepsilon$, and currently $\delta > \varepsilon$. Mathlib's 2 million lines of code represent an impressive achievement, but vast areas of mathematics remain unformalized.

This is where AI assistance becomes essential. Note that I'm now talking about a much more modest goal than proving new theorems: if AI can help with what's called ``autoformalization'' -- automatically converting already established natural language mathematics into formal proofs -- then $\varepsilon$ could at least get on par with $\delta$. And then, once Mathlib is sufficiently comprehensive, we might finally 
get $\varepsilon$ to actually exceed $\delta$, allowing one to
work directly in the formal system to discover new mathematics.

How might such AI-assisted formalization work in practice?

\section{A Path Forward: Quasi-Autoformalization.}\

The answer, I believe, may lie in what I call ``quasi-autoformalization'': a system of cooperating AI agents, orchestrated to translate natural language mathematics into formal proofs.
The key word here is ``quasi.'' For my purposes as a research mathematician (rather than an AI researcher), I don't require the system to work completely autonomously. I'm perfectly happy to intervene whenever I see a quick way to help prune the search tree and move the process along. My goal is not full automation but rather a collaborative and dramatically accelerated workflow. Here's an overview, in broad strokes, before getting to the details:

\begin{figure}[h]
    \centering
    \includegraphics[width=0.8\linewidth]{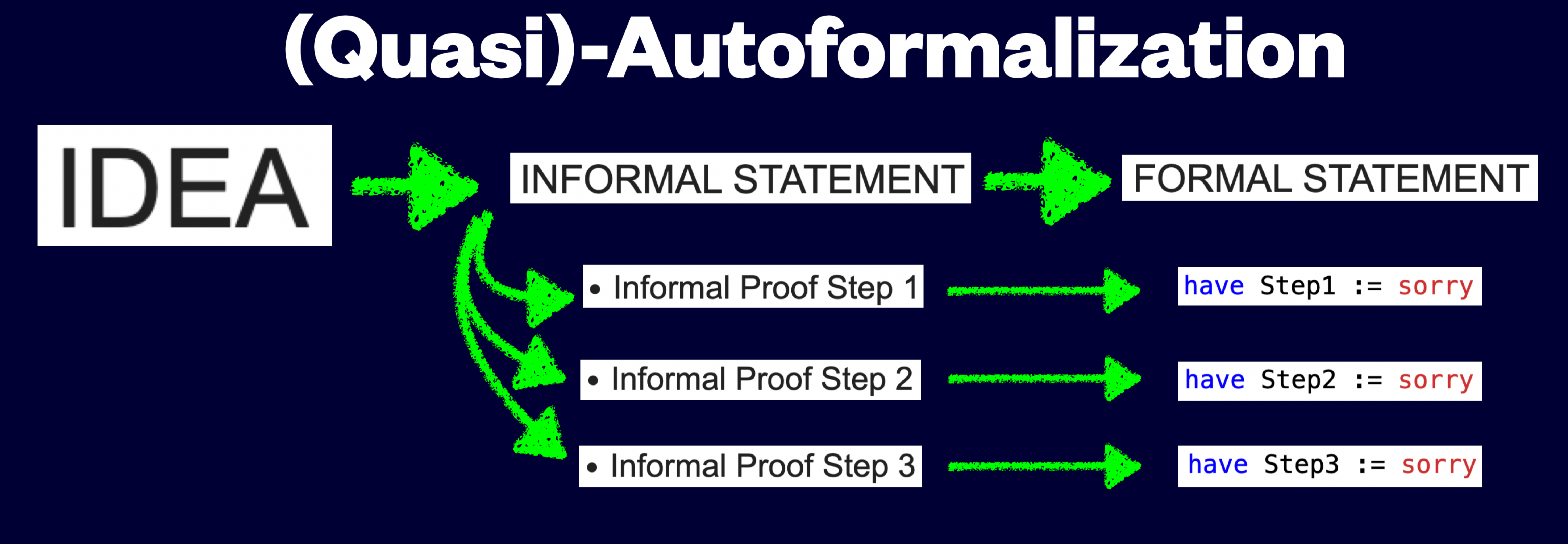}
    \caption{A schematic for Quasi-Autoformalization}
    \label{fig:quasi}
\end{figure}

Note that a variety of approaches to AI-assisted formalization have been explored in the past, from early work on automated theorem proving with human-readable output \cite{GanesalingamGowers2017} to recent LLM-based systems like Draft-Sketch-Prove \cite{DraftSketchProve2023} (focused on solving new problems from scratch) and tools like Lean Copilot \cite{LeanCopilot2024} (providing suggestions within the proof assistant). The quasi-autoformalization framework I propose differs in its explicit multi-agent decomposition with iterative refinement, designed specifically to accelerate Mathlib growth rather than (necessarily) achieve full autonomy.

\subsection{The Architecture: Four Agents.}
The system consists of four components working together:

\begin{itemize}
\item \textbf{Decomposer}: Takes mathematics at the level of ``Ideas'' and refines it into natural language statement and proof steps at an appropriate level of granularity.

\item \textbf{Translator}: Converts natural language statements and proof steps into formal Lean code, producing both the theorem statement and a ``scaffold'' of \lean{have} statements (intermediate claims) that are initially \lean{sorry}'d out (that is, their proofs are postponed without setting off Lean error messagees).

\item \textbf{Solver} (or ``Closer''): Attempts to prove each \lean{have} statement, replacing \lean{sorry}s with complete formal proofs.

\item \textbf{Conductor}: Orchestrates the entire process, identifying failures and deciding which components need to re-run with adjusted parameters.
\end{itemize}

Let's examine each component in detail.

\subsubsection{The Decomposer: From Ideas to Scaffold.}

Natural language mathematics, whether in undergraduate textbooks or arXiv papers, is really written at the level of ``Ideas.'' Even something as clear (to a trained mathematician) as a proof in Baby Rudin is nowhere near precise enough to formalize directly. 

\begin{figure}[h]
    \centering
    \includegraphics[width=0.5\linewidth]{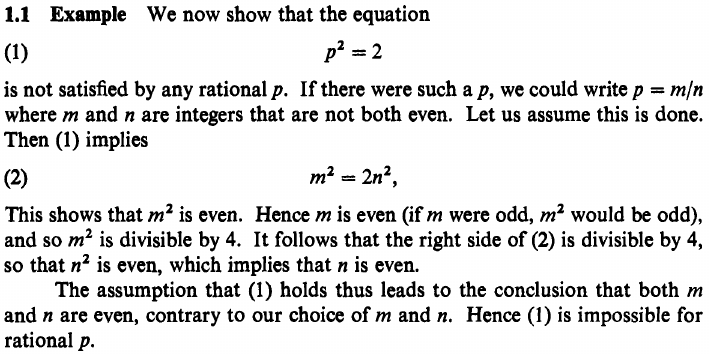}
    \caption{Taken from Rudin \cite[p.2]{Rudin1976}}
    \label{fig:Rudin}
\end{figure}

Formalizing this simple argument requires many more steps that are hidden in  natural language because they're too obvious to a human to state but are necessary for completeness, such as: why can we make it so that $m$ and $n$ are  not both even? Why if $m$ is not of the form $2k$, then it must be of the form $2k+1$? Etc.

\begin{figure}[ht]
    \centering
    \includegraphics[width=0.5\linewidth]{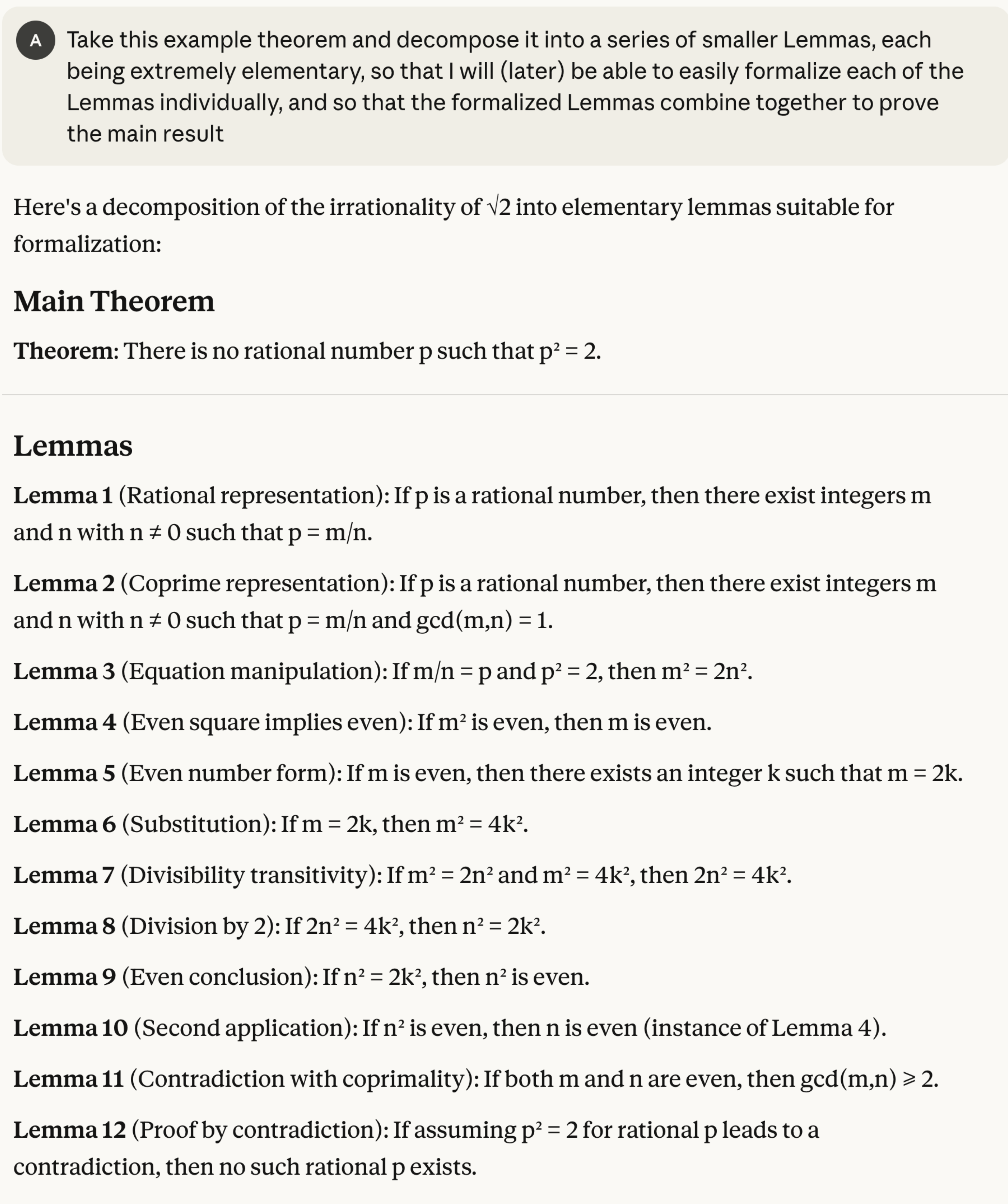}
    \caption{A decomposition of Rudin's argument using Claude AI.}
    \label{fig:decomp}
\end{figure}

The Decomposer's job is to break down mathematical ideas into sufficiently small, explicit steps that can be formalized. The granularity depends on what's already in Mathlib and on how sophisticated the Solver is. This is an iterative process: if the Solver fails, the Conductor may ask the Decomposer to break things down further. \figref{fig:decomp} shows a sample decomposition in this example, by Claude AI.

\subsubsection{The Translator.}

\begin{figure}[ht]
    \centering
    \includegraphics[width=0.45\linewidth]{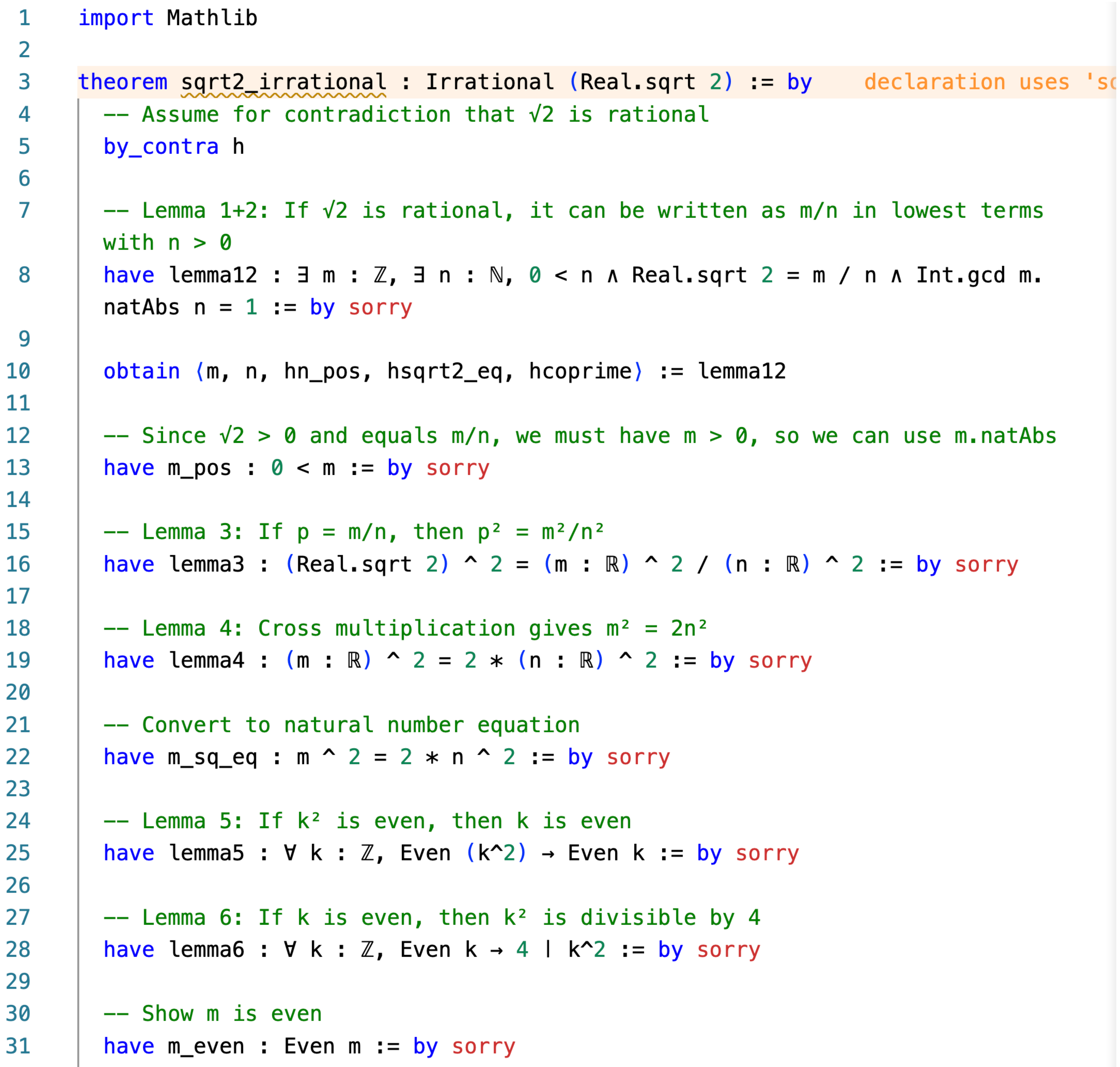}
    \includegraphics[width=0.45\linewidth]{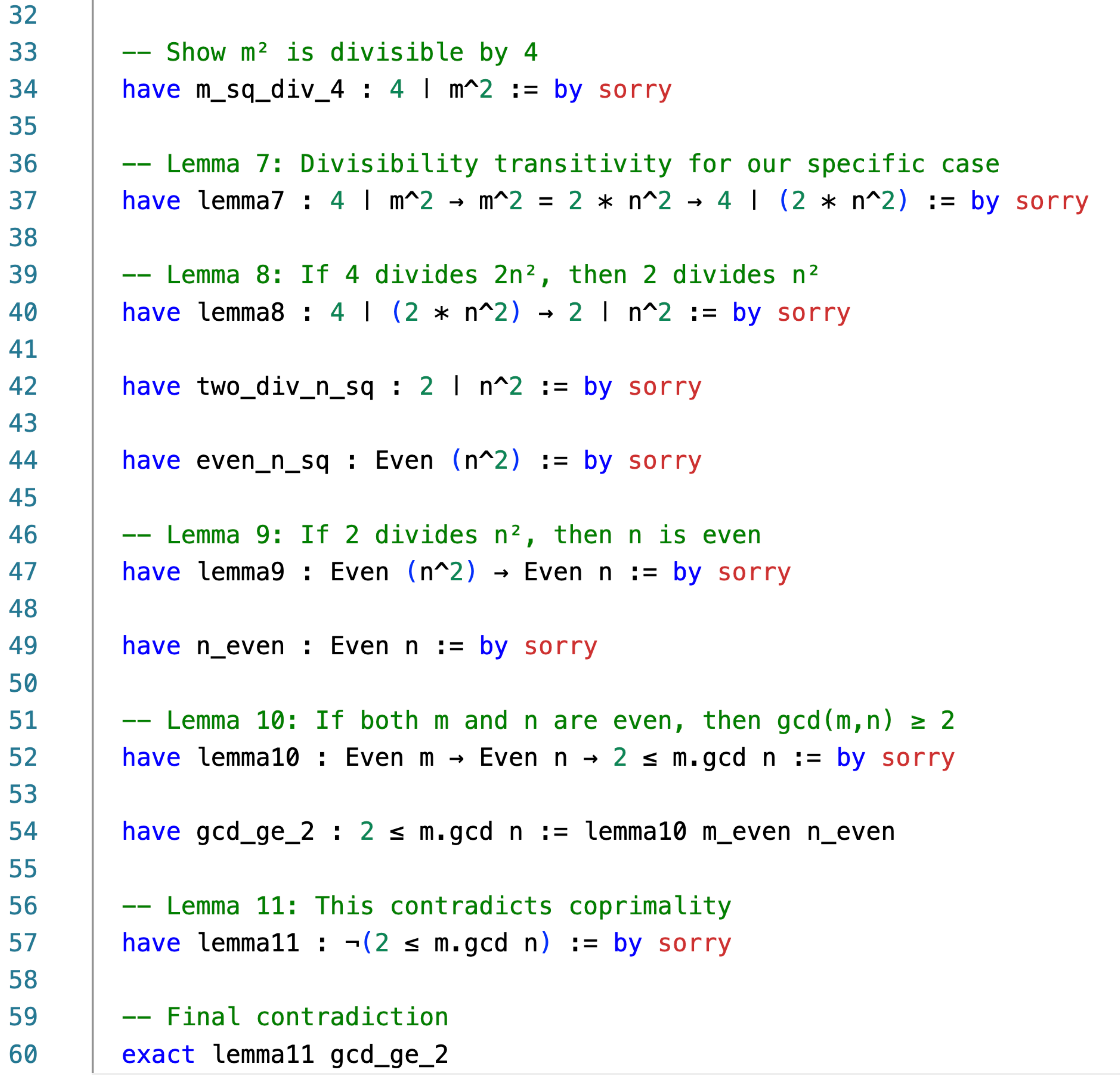}
    \caption{A translation of the natural language statement and lemmas into a scaffold for a formal proof}
    \label{fig:trans}
\end{figure}

The Translator's job is to take the Decomposer's small natural language lemmas, and convert each  into a small formal Lean \lean{have} statement, giving a formal scaffold of the main theorem. \figref{fig:trans} shows what it might look like in our running example, again using Claude AI. We postpone to \secref{sec:trans} a more detailed discussion of the manifold difficulties this particular agent will encounter.

\subsubsection{The Solver: Closing Goals.}

Once the Translator produces a scaffold with \lean{sorry}'d \lean{have} statements, the Solver attempts to clear each \lean{sorry}. This is where automated theorem proving systems demonstrate their power: they can search through vast spaces of possible proof strategies, chaining together lemmas from Mathlib, and applying tactics until goals are closed.

The Solver doesn't need to succeed on every goal. Some \lean{sorry}s may remain, requiring either human intervention, further decomposition by the Decomposer, or refinement of the statement by the Translator. The system's value lies in automatically closing the goals that \emph{can} be closed, dramatically reducing the human effort required.

\begin{figure}
    \centering
    \includegraphics[width=0.45\linewidth]{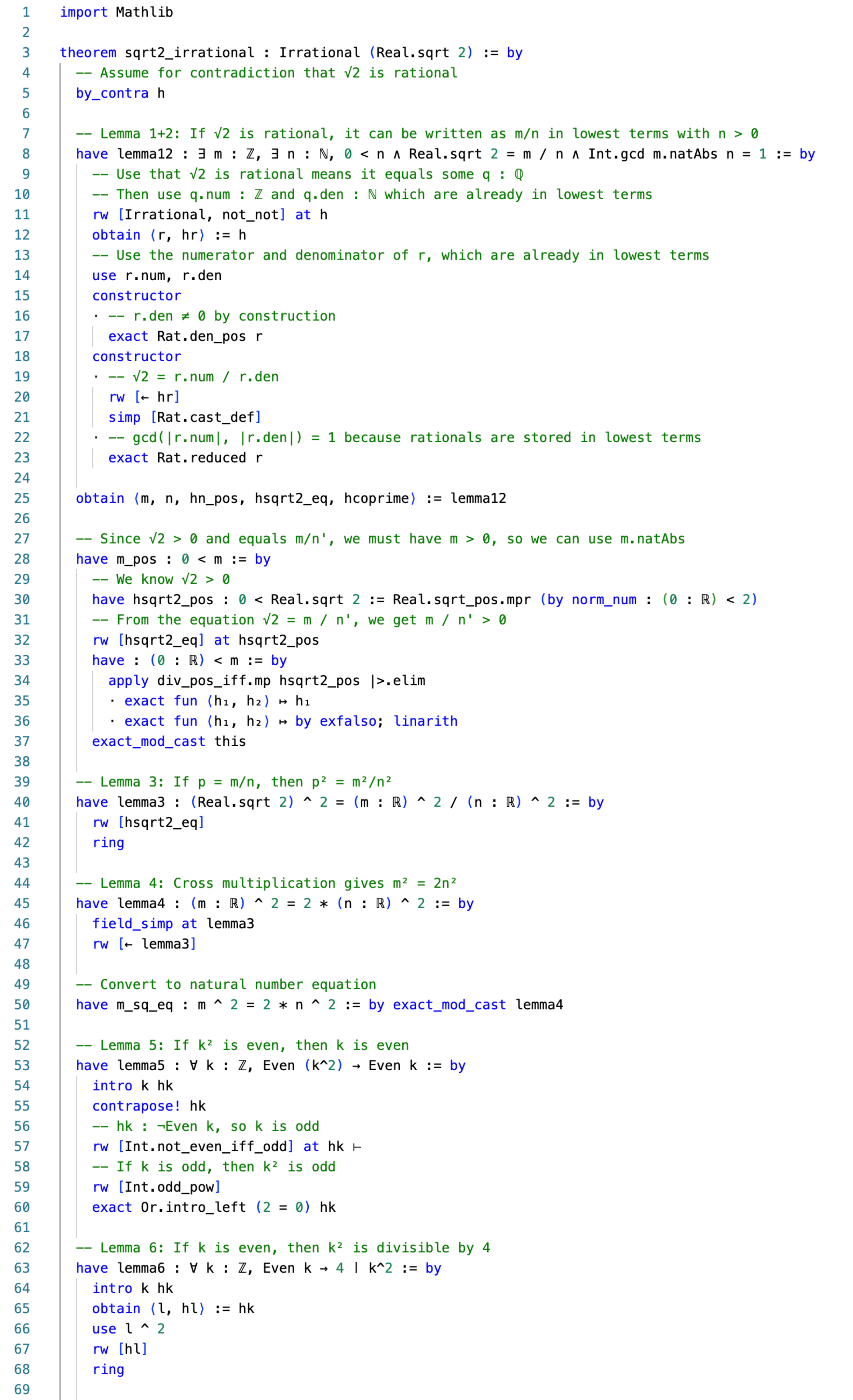}
    \includegraphics[width=0.45\linewidth]{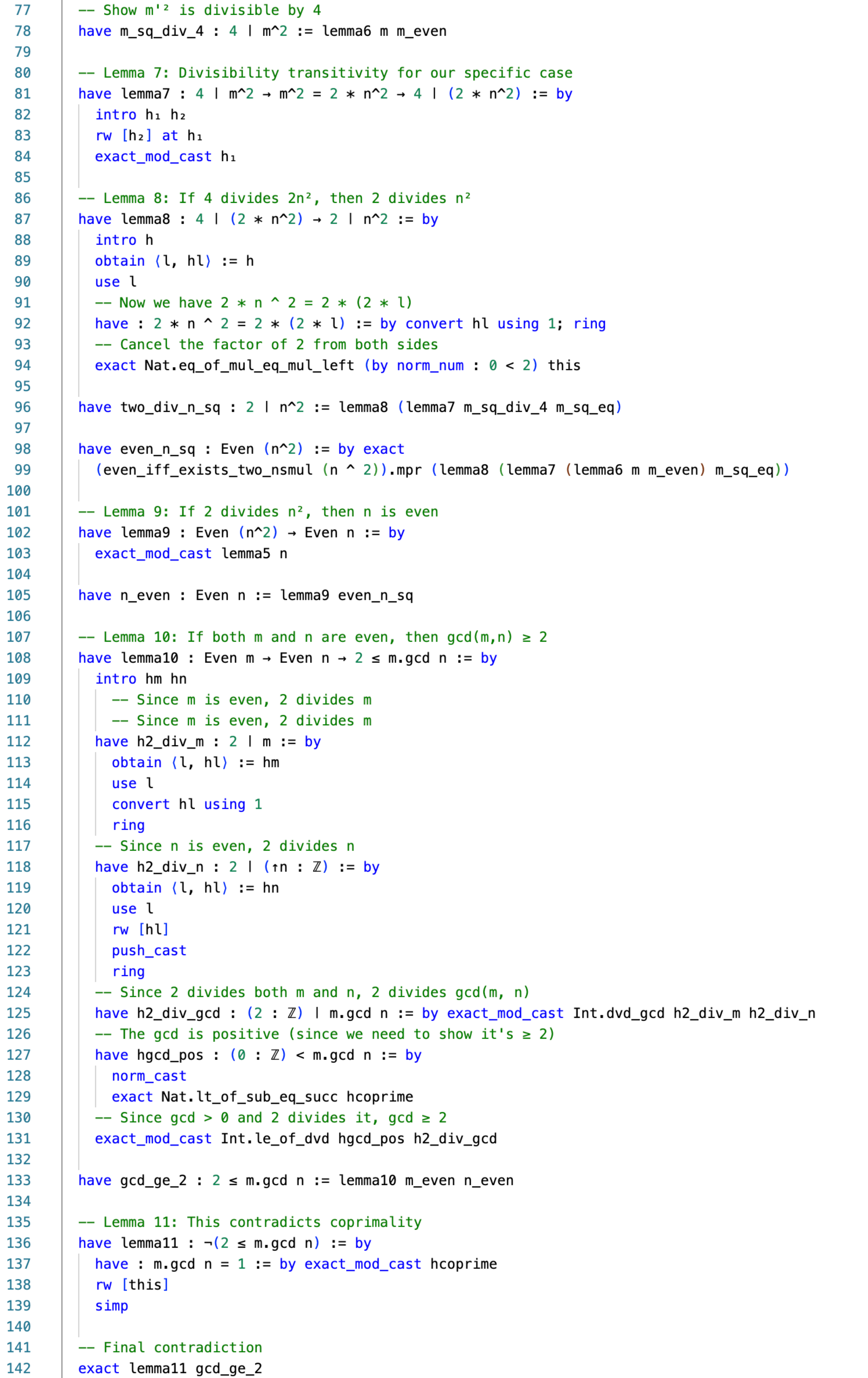}
    \caption{The Closer successfully solved all the remaining goals.}
    \label{fig:closer}
\end{figure}

Different solving strategies may be appropriate for different types of goals. Simple algebraic manipulations might yield to basic tactics, while more subtle arguments may require sophisticated search procedures or even reinforcement learning approaches trained on large corpora of formal proofs. \figref{fig:closer} shows a Closer applied to our example (with some help from this human).

\subsubsection{The Conductor: Orchestrating the Workflow.}

The Conductor is the system's coordinator. When formalization fails (as it inevitably will on first attempts), the Conductor must diagnose why:

\begin{itemize}
\item Is a goal beyond the Solver's current capabilities? Flag for further decomposition.
\item Was a statement mistranslated? Have the Translator try again with additional context.
\item Does this result perhaps already exist in Mathlib? Perform a Google search (via RAG -- \emph{retrieval augmented generation} -- or other methods) to check whether the theorem (or something closely related) is  already formalized in the latest version of the library.
\end{itemize}

The Conductor also manages computational resources, deciding when to re-run components versus requesting human input. In the ``quasi'' approach, a human can serve as super-Conductor, making these decisions based on mathematical insight.

\begin{figure}
    \centering
    \includegraphics[width=0.7\linewidth]{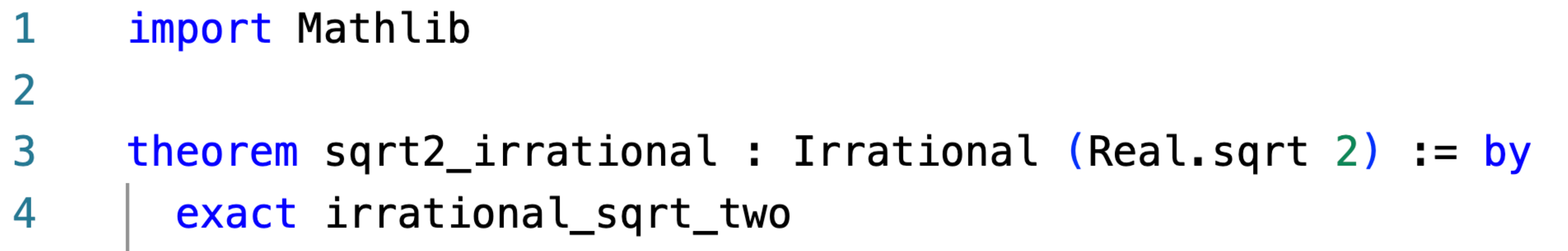}
    \caption{A search of the library should find this result.}
    \label{fig:conductor}
\end{figure}

In the case at hand (see \figref{fig:conductor}), the Conductor's RAG lookup of the library could have found a theorem called \lean{irrational_sqrt_two}, thereby reducing 140 lines of code to a single citation of Mathlib's much more elegant proof.

\subsection{Optimism...}

Does this workflow actually deliver? The first time I drew the schematic in \figref{fig:quasi}, as soon as it became clear in my mind, was in May 2025 at the Institute for Advanced Study (where I was spending my sabbatical) over lunch with DeepMind’s Thomas Hubert, a leader on the AlphaProof project, and Princeton's Sanjeev Arora and Chi Jin, leaders of the open source Goedel Prover \cite{GoedelProver2025, Goedel2Prover2025}. 
Both AlphaProof and Goedel Prover represent the current state of the art in automated theorem proving, with the latter competing in the Open Source category.
The day before, I had the opportunity to test AlphaProof with Thomas against some outstanding tasks in the ``PrimeNumberTheoremPlus Project'' (PNT+) which  Terry Tao and I co-organize \cite{PNTPlus2024}. The project's goal is to formalize significant results in analytic number theory, including the prime number theorem, Dirichlet's theorem on primes in progressions, and the Chebotarev density theorem. Seeing a real ``Closer'' in action for the first time made the potential of this approach suddenly concrete.

The Solver faces a genuinely difficult task: even goals that seem trivial to a human mathematician can require intricate formal reasoning. For example, I had scaffolded a lemma requiring a \lean{have} statement that the Riemann zeta function is at least 1 in absolute value in a neighborhood of $s=1$ (excluding the point itself). This should be ``immediate'' to a human: it follows from \lean{riemannZeta_residue_one}, the theorem already in Mathlib stating that $\zeta$ has a pole at $s=1$ (with residue $1$). But if proving that $\sqrt{2}$ is irrational took 140 lines of formal code, one can only imagine how much more complicated this ``trivial'' task becomes.

\begin{figure}[ht]
    \centering
    \includegraphics[width=0.8\linewidth]{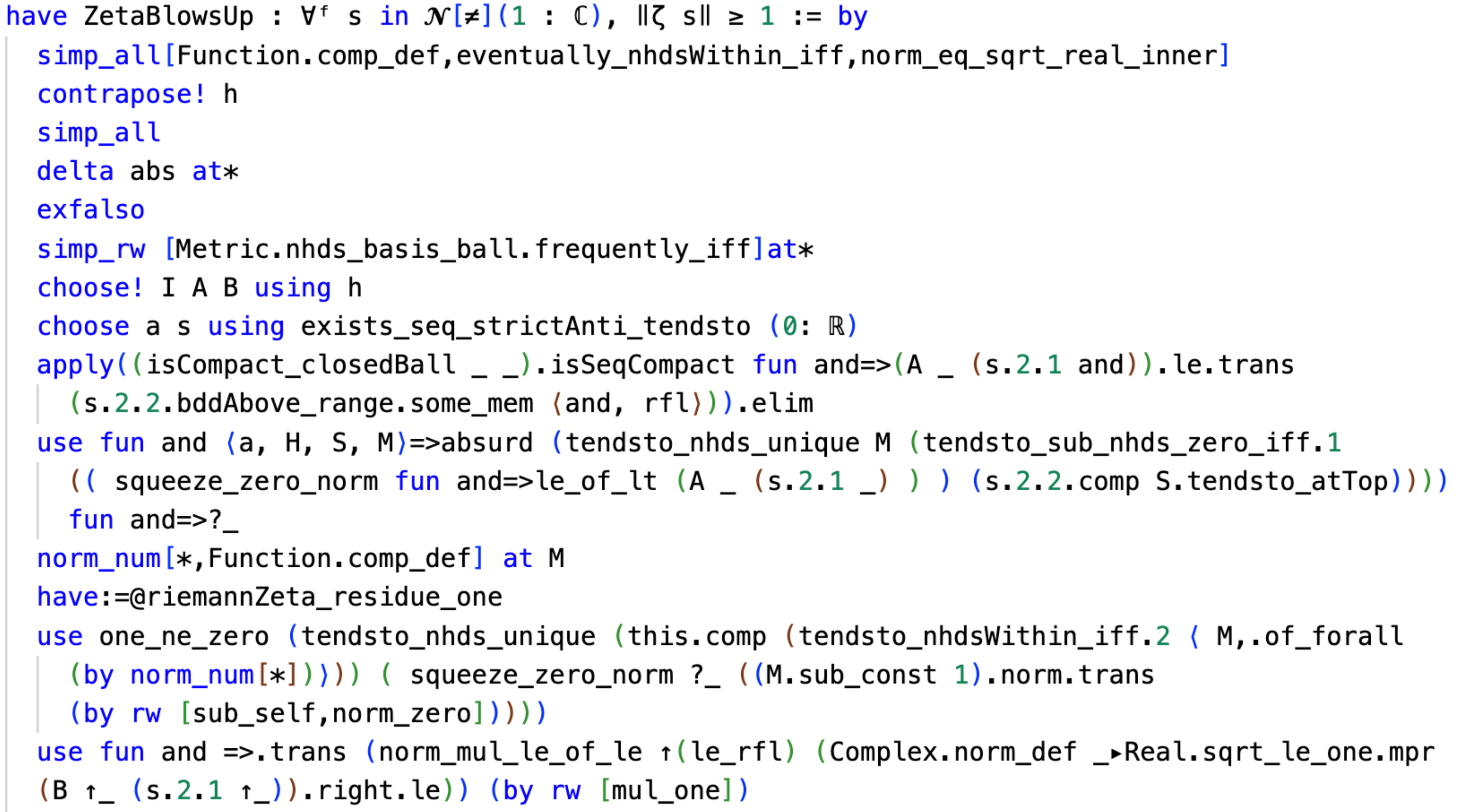}
    \caption{AlphaProof's solution to a \lean{have} statement, from \cite{PNTPlus2024}}
    \label{fig:zeta}
\end{figure}

And yet AlphaProof succeeded brilliantly! This despite being trained only on high school IMO problems (which never involve ``eventually filters,'' limits, or complex analysis), it chained together an extremely long sequence of formal steps that Lean accepted. The resulting proof was barely comprehensible, even to an expert human formalizer, see \figref{fig:zeta}.

This seems to reflect AlphaProof's training mechanism: adapting the AlphaZero reinforcement learning protocol, it ``plays'' mathematics against itself, discovering highly non-standard patterns for proving theorems through self-improvement.
Since my work wasn't meant to directly contribute to Mathlib (yet), I didn't care what the formal proof looked like, as long as the Lean code compiled.
 Since I knew the result followed trivially in natural language from already-proved theorems, Lean's acceptance was sufficient, and I could move on to the next goal.

In September 2025, shortly before this document's  completion, 
Morph AI (now Math, Inc.), led by Jesse Han, Jared Lichtman, and Christian Szegedy, announced another significant milestone: their ``Gauss'' agent autoformalized in Lean the classical error term for the Prime Number Theorem. It built on the foundations of PNT+, and made substantial progress formalizing complex analysis results such as the Borel-Carathéodory theorem and extending it to logarithmic derivatives of meromorphic functions. In their experiment, a mathematician (Jared) served as Decomposer while Gauss handled translation and solving autonomously -- a division of labor that may itself become fully automated in future iterations.

While these successes point to the framework's viability, significant challenges remain.

\subsection{... and Difficulties.}

Of course it is usually not so straightforward to get this workflow to succeed; let's highlight a number of particular difficulties.

\subsubsection{The Challenge of Closing Goals.}

As we saw with AlphaProof's success on the zeta function goal, closing formal goals is possible but requires sophisticated systems. 
%Two 
A particular engineering feature proved especially valuable here:

%First, AlphaProof could run in the background while I continued editing Lean code. This seems trivial, but it significantly improved my workflow. (The standard search tactics like \lean{exact?}, \lean{apply?}, \lean{rw?}, and \lean{simp?} don't work this way; you must wait until they finish, and touching a key restarts their search from scratch.)

Given any goal, AlphaProof simultaneously searches for a proof of the statement \emph{and its negation}. Their original motivation for engineering this was practical. Many IMO problems ask for answers, not just proofs. If you don't know the precise statement you're trying to prove, you cannot even begin working formally. 
For example, a problem might ask for the set of integers satisfying some property.
You might formalize this statement as $\exists S \subseteq \mathbb{Z}, S = \{n \in \mathbb{Z} \mid \text{property}(n)\}$. But this challenge is \emph{trivially} solvable: just set $S := \{n \in \mathbb{Z} \mid \text{property}(n)\}$, and by definition, it's equal to itself! To make a nontrivial statement requires first determining what $S$ is, say, $S = \{1, 2, 4\}$, and only then trying to prove that $\{1, 2, 4\} = \{n \in \mathbb{Z} \mid \text{property}(n)\}$.

AlphaProof's approach: send the natural language problem to an LLM like Gemini, requesting roughly 100 potential solutions ($S = \{1, 2, 3\}$, $S = \{2, 3, 4\}$, etc.). For each guess, AlphaProof attempts to prove or disprove the claim. In practice, it quickly disproves approximately 98 of them, then focuses on the remaining candidates.

This feature proved invaluable in our experiments for a completely different reason: it catches translation errors.

\subsubsection{The Subtleties of Translation.}\label{sec:trans}

The Translator component presents particularly subtle challenges. Getting formal statements to say \emph{exactly} what we intend to is sometimes as hard as proving them!
 In my experiments with AlphaProof on the PNT+ Project, embarrassingly many of the \lean{have} statements I scaffolded were resoundingly \emph{disproved} by the system. Each time, I could usually quickly identify what assumption I had omitted, adjust the statement, and iterate. 
 This back-and-forth between proposing statements and checking whether they can be established (or refuted) was essential to getting the formalization right.
 
For an illustration, consider formalizing the Riemann zeta function. Many think of it as $\zeta(s) = \sum_{n=1}^{\infty} 1/n^s$, but this is really the Euler (or Bernoulli) ``zeta function''. Riemann's insight was that the ``right'' way to define the zeta function was by first making a theta function, taking its Mellin transform, and dividing by Gamma:
\begin{align}\label{eq:zeta}
\zeta(s):=\frac{\pi^{s/2}}{\Gamma(s/2)}\left[\int_1^\infty\left(2\sum_{n=1}^\infty e^{-\pi n^2 u^2}\right)(u^s+u^{1-s})\frac{du}u-\frac1{1-s}-\frac1s\right]
\end{align}
This agrees, for $\Re(s)>1$, with the series, and hence represents its meromorphic continuation.

Now, in most Interactive Theorem Provers, including Lean, functions are \emph{total} (that is, defined on all inputs). So although you're not supposed to divide by zero, the expression $1/0$ must have some value (called a ``junk'' value, because you should never have to know what it is); it cannot be left undefined. In Mathlib, the junk  value assigned  to $1/0$ is $0$. This makes it possible to prove theorems like $(a+b)/c=a/c + b/c$ without requiring a proof that $c\ne0$; the junk value is such that the statement is true either way. The rationale for this is simple and very practical: If you required denominators to be nonzero before you allowed division, you would quickly be swamped with a cascade of trivial theorems needing to be proved again and again, for every arising division. Instead, the practice is to push the need to prove nonzeroness as late as possible; for example, the theorem that $a / b \times b = a$ \emph{does} (as it must) require as a hypothesis that $b\ne0$.

What does this have to do with the zeta function? David Loeffler and Michael Stoll formalized it and other $L$-functions for Mathlib \cite{LoefflerStoll2025}, and discovered the following very curious phenomenon.
Since functions are total, the zeta function maps all of $\C$ to $\C$, and hence has some junk value at its pole at $s=1$. What is the value?

While the question seems trivial, note that
if that junk value were zero, the Riemann Hypothesis, as is commonly stated would be trivially false! The standard statement says that zeros in the critical strip (say $\Re(s) \in [0,1]$, avoiding the trivial zeros at negative even integers) satisfy $\Re(s) = 1/2$. But if $\zeta(1) = 0$, we would have a counterexample. Imagine an AI announcing that it has disproved the Riemann Hypothesis with a petabyte-long incomprehensible proof whose essence is that $\zeta(1) = 0$.

Fortunately, Loeffler and Stoll worked out carefully how the junk value propagates through the construction \eqref{eq:zeta}, and found  that $\zeta(1) \neq 0$. So the standard statement of RH requires no adjustment. But this example illustrates just how  subtle and treacherous translation can be.

When I discussed these translation difficulties with Christian Szegedy -- a visionary scientist who was among the first to strongly advocate for autonomous formalization as a means of building large mathematical libraries \cite{Szegedy2020}, and who leads the team developing the Gauss agent mentioned above -- his response shocked me.

\subsubsection{Mathematics Is Robust.}

Christian told me, ``It doesn't matter if the formalized statements and definitions are wrong! Don't you \emph{believe} in mathematics?''

His point was profound: mathematics, in practice, is \emph{robust}. Definitions and theorems emerge through years of refinement, not arbitrary invention. We routinely get things wrong initially and refine the theory when necessary.

There are numerous examples in the history of mathematics that illustrate this point. In topology, we started by working with intervals on the real line (open and closed). Then we realized that this had nothing to do with the real line, and extended the idea to balls in arbitrary metric spaces. Eventually, Hausdorff showed that you don't need a metric at all! You simply declare which sets are open and require the right axioms, and can do abstract topology. Each iteration (``refactor'') improved our understanding of the notion of Openness.

Gauss's original statement of the Class Number One Problem wasn't quite right \cite{Goldfeld1985}. The Poincaré conjecture wasn't at first stated exactly as we now understand it. Mathematicians regularly get statements slightly wrong, even if the spirit of the idea is eventually validated. As long as we have mechanisms to refactor, we recover correct formulations \emph{in the long run}.

Christian's point was that the same will hold for formalization. Mathematics is self-correcting; wrong statements that are useful will be discovered and corrected through use! (And wrong statements that are useless won't be corrected, but they'll be harmless since no one is building on them.) 

For a concrete instance of this, imagine stating and proving the so-called Archimedean Property (that no matter how small a real number $\varepsilon>0$ may be, you can always find a natural number $N$ large enough that $1/N$ is even smaller) as follows:
\begin{lstlisting}
    theorem ArchimedeanProperty : ∀ (ε : ℝ), ε > 0 → ∃ (N : ℕ), 1 / N < ε 
\end{lstlisting}
And suppose that you  want at some later point to use this fact to prove that the sequence $a_n = 1/n$ converges to $0$:
\begin{lstlisting}
    def SeqLim (a : ℕ → ℝ) (L : ℝ) := ∀ ε > 0, ∃ N, ∀ n ≥ N, |a (n) - L| < ε

    theorem OneOverN (a : ℕ → ℝ) (hypothesis : ∀ n, a (n) = 1 / n) : SeqLim a 0
\end{lstlisting}
You will at some point need to leverage the fact that $n \ge N$ into the fact that $1/ n \le  1/N$. But this is not, in general, true: if $N=0$, then $1/N=0$, as we discussed already; so the latter inequality can fail! The culprit is the statement of \lean{ArchimedeanProperty}, which should insist on $N$ being strictly positive. In fact, the version stated above has the following vacuous proof:
\begin{lstlisting}
    theorem ArchimedeanProperty : ∀ (ε : ℝ), ε > 0 → ∃ (N : ℕ), 1 / N < ε := by
      intro ε hε
      use 0
      simp_all
\end{lstlisting}
That is, let $\varepsilon$  be given and assume (\lean{hε}) that $\varepsilon>0$. Set $N=0$; then $1/N<\varepsilon$ follows from the assumptions.

Again, it is these kinds of situations (and others in much less obvious contexts!) that make me feel that, as a mathematician, I must hold \emph{myself} personally responsible for any formal statements, whether translated by hand or by machine. Christian's counterpoint is intriguing: he claims that human mathematicians don't need to learn Lean at all; they can simply let autoformalizers do their thing, and trust that \emph{in the long run}, the important issues will sort themselves out.

However, recent examples of agentic AI behavior suggest caution about such trust. When I asked Claude to compute $0.99^{1000}$, it correctly recognized this as a deterministic question, wrote perfect Python code to perform the calculation, executed it internally, and reported the exact result. This demonstrates how tool use can bridge the stochastic-deterministic gap. But during the aforementioned DeepMind visit to IAS, we were told about a more troubling example: Gemini was asked to compute the genus of a particular curve. It produced perfect Sage code and reported that the genus was one. But the mathematicians asking  the question knew that the correct answer was two! Had they just discovered a bug in Sage? Upon further investigation, they found the following.
Their experimental version of Gemini was equipped with a safety feature to disable external tool access, and it was accidentally turned on; so Gemini \emph{never called} Sage at all! The model was being trained to report Sage's result but couldn't, so it decided to simply \textbf{pretend} to execute the Sage code and report its best guess at the result!

This illustrates why, while I value and respect Christian's perspective, I cannot endorse it myself: even when an AI claims to have used formal verification tools, we cannot blindly trust that it actually did so. Whether one adopts my cautious view or Christian's optimistic one may depend on temperament, but both viewpoints agree on the fundamental goal: accelerating the growth of formalized mathematics.

\section{The Adoption Question: When Will Formal Mathematics Win?}\

If the above pans out, then in the near future, Mathlib will grow large and robust enough that I can finally not only state and prove my theorems in Lean, but actually work on original research directly in the formal system, rather than first discovering results in natural language and then attempting to translate them. Will others follow suit?
The path to widespread adoption may follow a pattern we've seen before. Consider what happened with \LaTeX{}.

\subsection{The ``Knuth Factor''.}

In 1978, Don Knuth released \TeX{}. At that time, mathematicians wrote their papers by hand and gave them to secretaries for typesetting. After waiting weeks or months, they received something approximating their original text, then iterated a few times (if they cared to), correcting typos introduced in the process.

In the 1980s, Mike Spivak promoted AMS\TeX{}, yet still relatively few mathematicians adopted it. The system remained too cumbersome for all but the most dedicated adherents.
By the mid-to-late 1980s and early 1990s, \LaTeX{} emerged with beautiful, easy-to-use macros and automation. Almost everybody switched to self-typesetting. Eventually Overleaf arrived, eliminating (for many) the need to install software locally. You could do everything in a browser using a free web application.

But \LaTeX{} adoption wasn't driven solely by efficiency in producing final documents. It became an organizational tool for the research process itself. When working on a substantial result, mathematicians make incremental progress on various lemmas, subcases, and technical estimates. Keeping track of what has been established, how pieces interconnect, and which gaps remain becomes unwieldy when scattered across handwritten notes. \LaTeX{} allows you to maintain a living document that evolves with your understanding, where you can easily reorganize arguments, insert new lemmas in the logical sequence, and cross-reference results as the proof architecture develops. The act of typesetting becomes part of mathematical thinking, not merely its final presentation.

I define the ``Knuth factor'' as the ratio of time required to develop and document a mathematical result using \LaTeX{} (once you've learned the syntax with its dollar signs, backslashes, and curly brackets) to the time required to develop and document the same result by hand. Around 1990, the Knuth factor dropped below 1. Shortly thereafter, nearly everyone switched -- without any coercion -- simply because it was the obvious way to speed up their workflow.

I expect the same will happen with Lean. Here's how.

\subsection{The Formalization Factor.}

In formalization literature, there's the so-called ``de Bruijn factor,'' which measures the ratio of the number of lines of formal code to lines of natural language proof \cite{deBruijn1980, Wiedijk2000}. Estimates vary, but factors of 4-10 are commonly cited. However, this metric, I believe, misses the point. 

Lines of code are no longer a meaningful proxy for human effort. LLMs can generate thousands of lines of code quickly, and automated solvers can fill in proofs that would have taken humans hours. What matters is not how many \emph{lines} are produced, but the amount of \emph{time} it takes mathematicians to do their work.

The proposed alternative metric is what I'll call the ``formalization factor'': the ratio of time required to discover and formalize a mathematical result (from initial idea through verified formal proof) to the time required to discover and write up the same result in natural language mathematics, typeset in \LaTeX{}. Crucially, this metric acknowledges that formal systems may accelerate not just the verification phase, but the discovery process itself. 
Just as \LaTeX{} became a tool for organizing evolving mathematical arguments, formal systems could provide structure and error-checking that accelerates the incremental development of complex proofs while making them more reliable.

This factor currently sits well above 1 -- perhaps 10, perhaps 100, depending on the subfield and the mathematician's familiarity with Lean. The difficulty isn't just writing proofs; it's that vast areas of mathematics aren't yet formalized enough to even state your theorems. You must either build the foundations yourself (prohibitively expensive) or wait for Mathlib to grow.

But once that factor drops below 1, I expect nearly everyone will voluntarily switch to working formally, just as they did with \LaTeX{}; no coercion will be necessary. It will simply be the obvious thing to do to speed up your workflow and increase confidence in your results.

\subsection{Conditions for Adoption.}

For the formalization factor to drop below 1, several conditions must be met:

First, Mathlib must be comprehensive enough that stating your theorem requires minimal foundational work. This is where quasi-autoformalization becomes essential: AI assistance can accelerate library growth so that $\varepsilon > \delta$, ensuring formalized mathematics keeps pace with (or exceeds) natural language mathematics. But library growth requires more than just correct proofs -- it demands polished, maintainable code. The AI must learn to write reusable mathematics that integrates cleanly with existing infrastructure, not merely produce working but unmaintainable formal arguments.

Second, the tools must become more accessible. Installation should be effortless. Imagine a browser-based environment like Overleaf, but for Lean, where you simply open a page and start working. (This nearly exists with platforms like \url{live.lean-lang.org} and GitHub Codespaces, though significant friction remains.)

Third, AI assistance must handle the tedious parts. Natural language interfaces could help write formal code; automated solvers could close routine goals; translation tools could convert mathematical exposition to formal scaffolds. The mathematician focuses on the creative work, identifying the key ideas and structuring the argument, while AI handles the mechanical formalization.

Fourth, and perhaps most importantly, the verification benefit must become immediate and tangible. Currently, formalization is an investment whose payoff is distant: you formalize now so that others can build on your work later, or so that you're certain your proof is correct. But if working formally means catching errors early, avoiding tedious case-checking, and confidently building on others' results, then the workflow itself becomes more efficient. The formalization factor drops not just because formalization gets faster, but because the entire research process improves.

When will this factor drop below 1? That depends primarily on the first condition: comprehensive libraries. With sustained AI-assisted formalization, Mathlib could reach the necessary scale within 5-10 years for many core areas of mathematics. For certain highly specialized fields, it may take longer. But once your field crosses the threshold, the switch could be rapid.

\section{Teaching.}\label{sec:Teaching}\

The adoption timeline I've outlined focuses on research mathematics, but teaching represents a parallel and potentially accelerating force. 
Among the most rapid converts to formalization are young people. To see why, imagine a world in which we tried to teach newcomers to chess purely by \emph{describing} moves in sequence, say using algebraic notation. The teacher says:

\begin{quote}
Ok class, the game went like this: 1.\ Nf3 Nf6 2.\ c4 g6 3.\ Nc3 Bg7 4.\ d4 0-0 5.\ Bf4 d5 (this is a transposed Grünfeld Defence) 6.\ Qb3 dxc4 7.\ Qxc4 c6 8.\ e4 Nbd7 9.\ Rd1 Nb6 10.\ Qc5 Bg4 11.\ Bg5 Na4!! Wow, can you believe that move?!
\end{quote}

Chess aficionados can effortlessly translate these symbols into an evolving series of positions, updating the board in their heads at every move. But for the rest of us, it’s vastly easier if you just \emph{show me the board}!

And yet this imagined world is \emph{precisely} how we currently teach almost all newcomers to theorem proving: entirely in natural language. At every step of a mathematical proof, the ``game board'' -- meaning the current objects, assumptions, and goal -- is changing. Professional mathematicians track this effortlessly (System I), but for beginners it requires deliberate, error-prone effort (System II). Think back to our Baby Rudin example: when proving that $\sqrt{2}$ is irrational, the underlying game board is shifting at each line of reasoning.

Of course, we would never dream of writing out the full board position after every single move in a natural-language proof. First, it would be unbearably cumbersome. Second, most of us got along ``just fine'' without it. In fact, one might argue that developing the ability to construct these invisible boards in one’s head is \emph{essential} to becoming a mathematician, just as a strong chess player must learn to visualize multiple positions in their mind’s eye.

But here is the crucial point: when teaching with Lean, the game board is \emph{always there}, automatically generated and continuously updated. Beginners don’t have to struggle to reconstruct it from terse prose -- they can simply look.

This is why, in my diagram in \figref{fig:Lean}, I drew an arrow from Mathlib to Teaching, but also a separate arrow directly from Lean itself. The Mathlib arrow is clear enough: teaching is easier when libraries are large and convenient. But the Lean arrow reflects something different: in teaching, one often does \emph{not} want the most general and abstract definitions as they appear in Mathlib. Just as we do not introduce high school students to Lebesgue integration before they have seen Riemann integration, so too in Lean it is often more effective to write out simpler, more concrete versions of definitions directly. For example, newcomers to Real Analysis should first encounter limits via the familiar $\varepsilon$--$\delta$ definition, rather than Mathlib’s general \lean{Filter.Tendsto}. In this way, Lean provides a flexible framework for pedagogy, allowing instructors to tailor the level of abstraction to the needs of their students.

\begin{figure}[ht]
    \centering
    \includegraphics[width=0.8\linewidth]{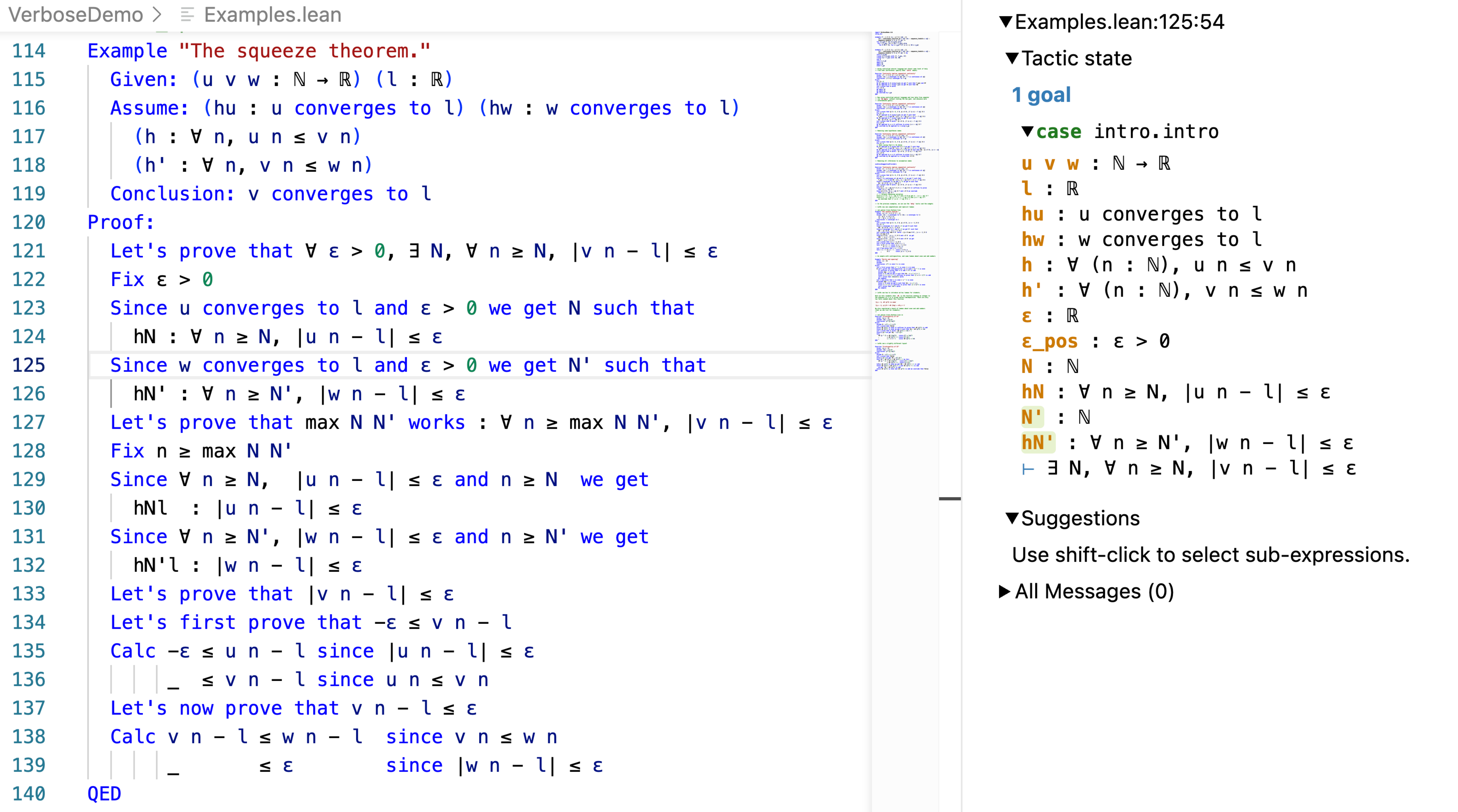}
    \caption{Massot's Verbose Lean wrapper for teaching Real Analysis in controlled natural language.}
    \label{fig:verbose}
\end{figure}

\begin{figure}[ht]
    \centering
    \includegraphics[width=0.8\linewidth]{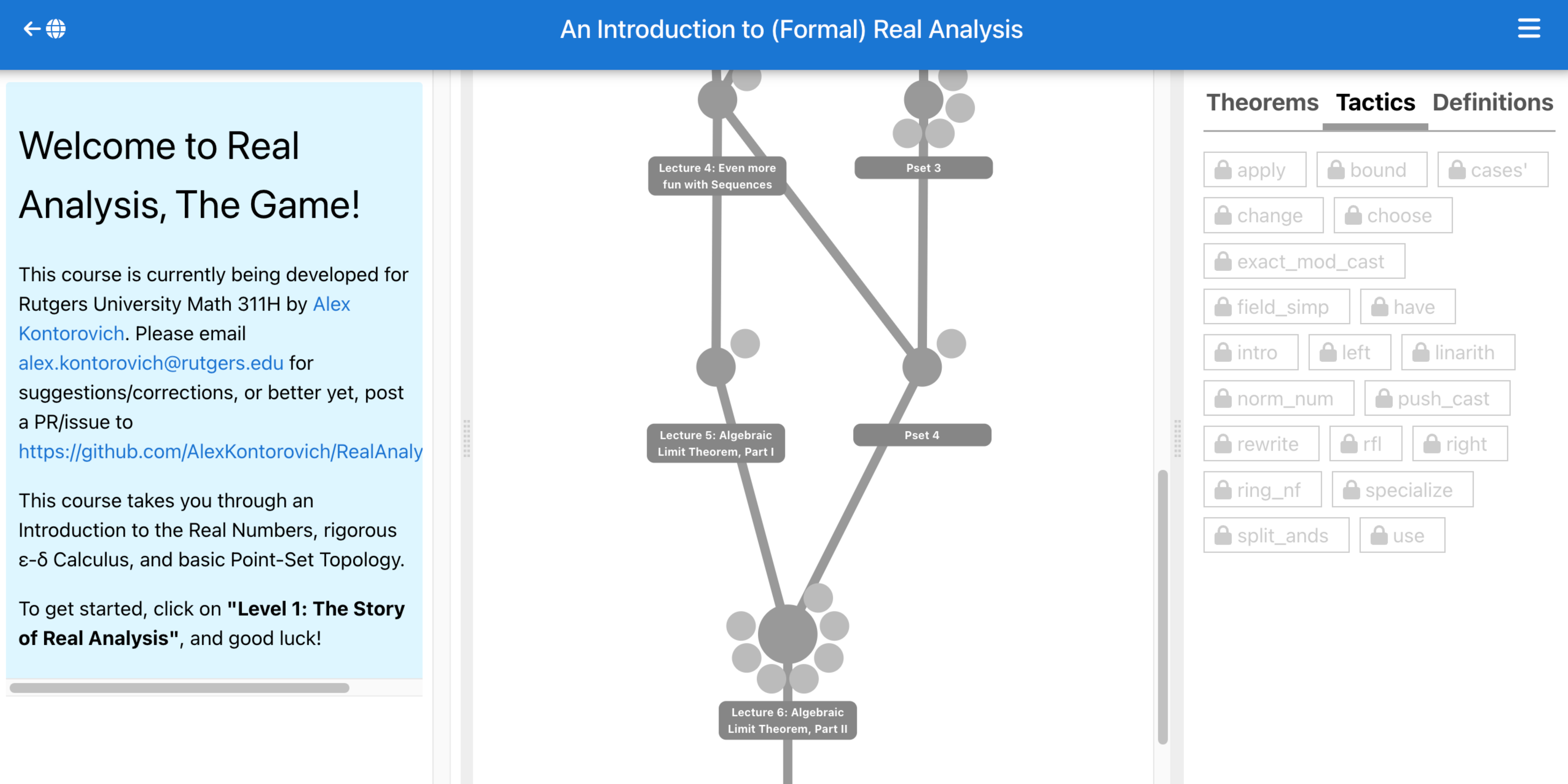}
    \caption{The Real Analysis Game, {\url{https://adam.math.hhu.de/\#/g/AlexKontorovich/RealAnalysisGame}}}
    \label{fig:game}
\end{figure}

This is an area ripe for experimentation. For but one example of very many, Patrick Massot \cite{Massot2024} developed a Lean wrapper called \textbf{Verbose}, which presents proofs in a controlled natural language, easing students into formal syntax, see \figref{fig:verbose}.
For my part, I am currently running an experiment in which I teach an undergraduate real analysis course as a kind of \emph{video game}, in the spirit of Buzzard's popular \emph{Natural Number Game}, see \figref{fig:game}.

These pedagogical experiments reflect a broader historical pattern. Before 1600, mathematicians were doing serious algebra, but it looked like this: 
\begin{quote}
You have some unknown quantity. When you create a second copy of this quantity and combine it with a known amount, then multiply this total by your original unknown quantity, you get a specified result.    
\end{quote}
Compare to the modern version:
$$x(x+a)=b.$$ 
The latter is, with some training, so much easier and immediate to understand.

Similarly, before formalization, we rarely named our hypotheses explicitly in proofs. I suspect this may change as formal systems influence mathematical exposition. Just as algebraic notation revolutionized how we think about equations, the precision required by formal systems may transform how we communicate mathematical arguments, leading us to speak more clearly and precisely even in natural language mathematics. The early signs suggest that Lean is not just a tool for research mathematics: it is poised to transform both mathematical pedagogy and mathematical discourse itself.

\section{Understanding and Communication.}\label{sec:Understanding}\

The adoption timeline and pedagogical benefits I have outlined focus primarily on efficiency and accessibility. But formalization may also fundamentally transform how mathematicians understand and communicate mathematical ideas themselves -- changes that could accelerate adoption for reasons that transcend mere productivity gains.

As Thurston eloquently argued in ``On Proof and Progress in Mathematics'' \cite{Thurston1994}, a formal proof (even in natural language) represents but the starting point of genuine mathematical understanding, not its culmination. The question is whether formal systems might actually enhance rather than constrain the development of mathematical insight and its communication.

\begin{figure}[h]
    \centering
    \includegraphics[width=0.8\linewidth]{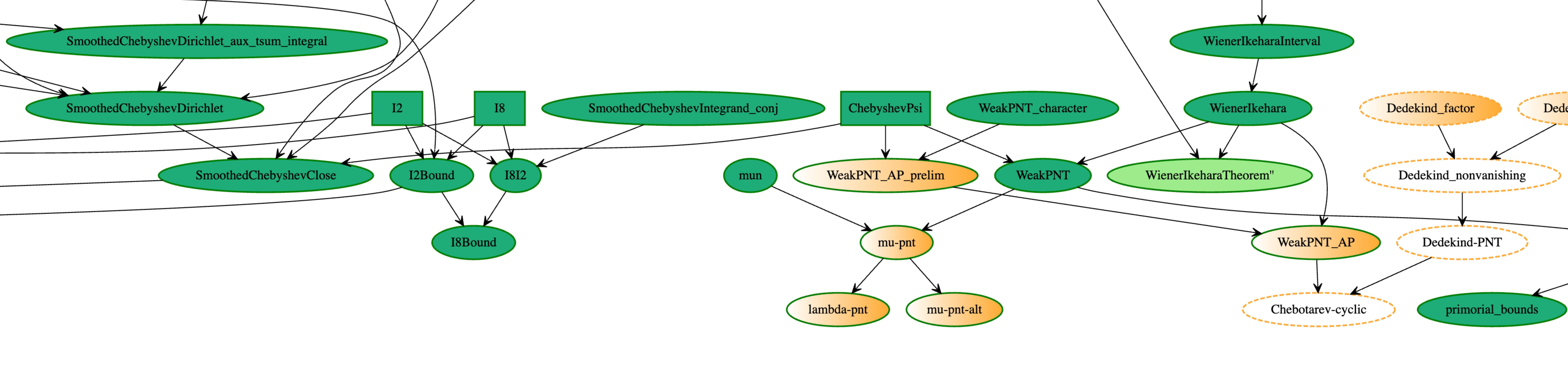}
    \caption{A portion of the dependency graph of \cite{PNTPlus2024}}
    \label{fig:dependency}
\end{figure}

Current large-scale formalization projects offer early glimpses of this potential. Massot's ``Lean blueprint'' technology seamlessly transitions between natural language \LaTeX{} and corresponding formal Lean code, while generating interactive dependency graphs that visualize the logical structure of complex arguments, see \figref{fig:dependency}. These graphs serve an immediate engineering purpose -- showing collaborators where ``leaves'' of the formalized tree require attention -- but they also point toward something more profound.

Consider how mathematical understanding might evolve when every theorem exists within an explicit, navigable web of logical dependencies. Rather than encountering results as isolated facts in linear exposition, mathematicians could explore arguments at variable levels of detail, zooming from high-level intuition down to granular formal steps as needed. A student learning algebraic topology could begin with the broad conceptual landscape, then drill down into specific homology calculations only when ready. An expert reviewing a paper could efficiently verify the overall logical structure, examining detailed proofs only for novel or suspicious claims.

This represents a qualitative shift from current mathematical communication, which remains fundamentally linear and static. Journal articles\footnote{For more on what the journal refereeing process may look like in the age of formalization, see \cite{Kontorovich2022}.} present arguments in fixed sequences, requiring readers to either accept claims on faith or laboriously verify every step. Formal systems could enable what we might call ``proof at multiple resolutions'' -- mathematical arguments that adapt to the reader's expertise and interests.

The implications extend beyond individual comprehension to collaborative mathematical practice. When disagreements arise about complex arguments, formal systems provide objective arbitration. When building on others' work, mathematicians could confidently incorporate results without fear of hidden errors. The cumulative nature of mathematical knowledge -- each generation building on previous insights -- could accelerate as the reliability of foundations increases.

This could enable a profound shift in mathematical responsibility itself. Currently, mathematical accountability is intensely personal and vertical. In graduate school, I was permitted to quote one result -- Deligne's proof of the Weil conjectures -- and expected to know everything else from the ground up. Even now, when I cite a theorem from another paper, I feel personally responsible for understanding its proof. We are like engineers who must mine their own metal, forge their own parts, and assemble everything personally.

But formalization, with its explicit specifications of every assumption and conclusion, could enable genuine modular mathematics. Imagine what theorems we could prove if we could confidently use off-the-shelf components! When formal statements specify precisely what they assume and what they guarantee, mathematicians could build on others' work with the same confidence engineers have in standardized components. This transition from vertical responsibility to horizontal collaboration could dramatically accelerate mathematical progress.

Yet significant challenges remain. Current formal languages prioritize logical correctness over human intuition, often obscuring the conceptual insights that drive mathematical progress. The risk is that formalization might fragment mathematical culture, creating a divide between those who work in formal systems and those who operate in natural language.

Whether formal systems ultimately enhance or constrain mathematical understanding will likely depend on our success in developing tools that preserve and amplify human mathematical intuition rather than replacing it. The experiments now underway in interactive visualization, natural language interfaces, and collaborative formal environments represent early steps toward what we might call the ``Paper of the Future'' -- mathematical communication that is simultaneously rigorous, accessible, and intellectually inspiring.

The stakes of this transition extend beyond mathematics itself. If formal systems can genuinely enhance mathematical understanding and communication, they may offer a model for rigorous reasoning in other domains where precision and cumulative knowledge matter. The shape of mathematics to come may prefigure broader transformations in how humanity approaches complex reasoning about the world.

\section{Conclusion.}\

The vision I have outlined -- of quasi-autoformalization accelerating library growth, formalization factors dropping below unity, and mathematical practice migrating toward formal systems -- represents both tremendous opportunity and considerable uncertainty. Whether this transformation unfolds as predicted depends on resolving fundamental tensions between the stochastic nature of current AI systems and the deterministic requirements of mathematical proof.

The evidence points in conflicting directions. On one hand, the trajectory from zero IMO capability in 2023 to gold medal performance in 2025 suggests that AI capabilities in mathematics are advancing with remarkable speed. The success of systems like AlphaProof in closing complex formal goals, even those requiring hundreds of intricate steps, demonstrates that sophisticated mathematical reasoning is within reach of current techniques.

Yet the core challenge remains: even a hypothetical AI system with 99.99\% accuracy per step would produce unreliable results when chaining thousands of reasoning steps together. The stochastic nature of large language models -- their fundamental reliance on probability distributions over next tokens -- appears deeply mismatched to the deterministic correctness that mathematics demands.

This returns us to Tim Gowers' prescient vision from 2000. Writing at the dawn of the new millennium, he predicted a coming ``Golden Age'' of human-computer collaboration in mathematics, but warned that such an era might be brief \cite{Gowers2000}:

\begin{quote}
The next stage might be one where only a very few outstanding mathematicians could discover proofs that were inaccessible to computers... In the end, the work of the mathematician would be simply to learn how to use theorem-proving machines effectively and to find interesting applications for them. This would be a valuable skill, but it would hardly be pure mathematics as we know it today.
\end{quote}

I offer a more optimistic perspective on such a Golden Age, should it indeed arrive. Mathematics has a rich history of problems that seemed intractably difficult but later proved more accessible than expected. The finite field Kakeya conjecture exemplifies this pattern: after years of incremental progress by leading experts, Zeev Dvir solved the full problem in 2008 \cite{Dvir2009} using elementary techniques from algebraic geometry, revealing that what seemed deeply difficult was actually rather straightforward.

This suggests that AI assistance may help us discover many open problems that are similarly more tractable than they appear, a kind of hidden low-hanging fruit throughout mathematics. The combination of AI's ability to explore vast solution spaces with formal verification's guarantee of correctness could unlock results that have remained elusive not due to fundamental difficulty, but simply because no human had the stamina to try the right combination of techniques from disparate fields.

Of course, what will likely endure after such successes are the genuinely hard problems; and here I fully expect challenges like the Riemann Hypothesis to remain formidable even for sophisticated AI systems. But if we do enter a period where formal, AI-assisted mathematics becomes the norm, the classical role of the ``heroic mathematician'' may persist in new forms: identifying which problems to pursue, recognizing when seemingly disparate areas might connect, and providing the conceptual insights that guide automated search.

The shape of mathematics to come will ultimately depend on choices we make today about how to develop these tools. If we prioritize pure automation over human insight, we risk creating systems that solve problems without advancing understanding. But if we succeed in building AI that amplifies rather than replaces mathematical intuition -- systems that handle the mechanical aspects of formalization while preserving space for human creativity -- we may witness not the end of pure mathematics, but its transformation into something more powerful and more beautiful than what came before.

Whether this vision proves prescient or merely optimistic, the mathematical community stands at a remarkable inflection point. The next decade will likely determine whether formal verification becomes as fundamental to mathematical practice as \LaTeX{} is today, or remains a specialized tool for the most dedicated practitioners. Either way, we are living through a pivotal moment in the history of mathematical reasoning, and the final chapter of this story remains to be written.

\section*{Acknowledgments.}
My path into Lean began with Kevin Buzzard’s online lectures and encouragement, and it continued only thanks to Heather Macbeth’s patient guidance. Paul Nelson was the first to show me the remarkable potential of large language models, and Christian Szegedy opened my eyes to the central importance of autoformalization. Along the way I have benefited from stimulating conversations with Sanjeev Arora, Jeremy Avigad, Tim Gowers, Thomas Hubert, 
Ian Jauslin,
Chi Jin, 
JJ Jung, Jeremy Kahn,
Patrick Massot, 
Konstantin Mischaikow,
Terry Tao, Akshay Venkatesh, and the vibrant Mathlib community. Above all, I am grateful to Leo de Moura, whose creation of Lean -- building on the rich tradition of type theory and proof assistants -- has opened formal mathematics to practicing researchers in ways previously unimaginable.

% SIAM recommends using BibTeX
% if using BibTeX
\bibliographystyle{siamplain}
\bibliography{AKbibliog}
\end{document}